\documentclass[11pt]{article}
\usepackage{amstext,amssymb,amsmath,amsbsy}

\usepackage{hyperref}
\usepackage{amscd}
\usepackage{amsfonts}
\usepackage{indentfirst}
\usepackage{verbatim}
\usepackage{amsmath}
\usepackage{amsthm}
\usepackage{enumerate}
\usepackage{graphicx}
\usepackage{subfig}
\usepackage{color}
\usepackage[OT1]{fontenc}
\usepackage[latin1]{inputenc}
\usepackage[english]{babel}
\usepackage{amssymb}
\newtheorem{theorem}{Theorem}
\newtheorem{lemma}{Lemma}

\newtheorem{proposition}{Proposition}

\newtheorem{remark}{Remark}

\newtheorem{definition}{Definition}
\numberwithin{equation}{section}
\numberwithin{theorem}{section}
\numberwithin{lemma}{section}
\numberwithin{notation}{section}
\numberwithin{proposition}{section}
\numberwithin{corollary}{section}
\numberwithin{corollary}{section}
\numberwithin{example}{section}
\numberwithin{definition}{section}
\numberwithin{remark}{section}

\newcommand{\proofend}{\hfill $\Box$ }
\newcommand{\mint}{\rule[1.1mm]{2.45mm}{.1mm}\hspace{-3.35mm}\int}

\newcommand{\dsp}{\displaystyle}

\newcommand{\supp}{\operatorname{supp}}

\newcommand{\dive}{\operatorname{div}}

\newcommand{\eps}{\varepsilon}

\newcommand{\loc}{_{loc}}

\newcommand{\mC}{\mathbb{C}}

\newcommand{\mN}{\mathbb{N}}

\newcommand{\mR}{\mathbb{R}}

\title{Cloaking via anomalous localized resonance for doubly complementary media in the quasistatic regime}

\author{Hoai-Minh Nguyen \footnote{EPFL SB MATHAA CAMA, Station 8,  CH-1015 Lausanne, hoai-minh.nguyen@epfl.ch}}

\begin{document}

\date{}
\maketitle

\begin{abstract} This paper is devoted to the study of  cloaking via anomalous localized resonance (CALR) in the two and three dimensional  quasistatic regimes. CALR associated with negative index materials  was discovered by Milton and Nicorovici in \cite{MiltonNicorovici} for constant plasmonic structures in the two dimensional quasistatic regime. 
Two key features  of this phenomenon are the localized resonance, i.e., the fields blow up in some regions and remain bounded in some others, and the connection between the localized resonance and the blow up of the power of the fields,  as the loss goes to 0. An important class of negative index materials for which the localized resonance might appear is the class of reflecting complementary media introduced in \cite{Ng-Complementary}.  It was showed in \cite{MinhLoc} that complementary property of media is not enough to ensure a connection between the blow up of the power and the localized resonance. In this paper, we study CALR for  a subclass of complementary media called  the class of doubly complementary media.  This class is rich enough to allow us to cloak  an {\bf arbitrary source} concentrating on an {\bf arbitrary smooth bounded manifold of codimension 1} placed in an {\bf arbitrary medium}  via anomalous localized resonance; the cloak is independent of the source. The following three properties are established for doubly complementary media: P1) CALR appears if and only if the power blows up;  P2) The power blows up if the source is  located ``near" the plasmonic structure; 
P3) The power remains bounded if the source is far away from the plasmonic structure. 
{\bf Property P2)}, the blow up of the power,   is in fact established for {\bf reflecting complementary media}. 
The proofs of these results are based on several new observations and  ideas. 
One of the difficulties in the study of this problem is to handle the localized resonance. To this end,  we extend  the reflecting and the removing localized singularity techniques introduced  in \cite{Ng-Complementary, Ng-Negative-Cloaking, Ng-superlensing}, and  implement the separation of variables for Cauchy problems for a general shell. The results in this paper are inspired by and imply recent ones of Ammari et al. in \cite{AmmariCiraoloKangLeeMilton} and Kohn et al. in \cite{KohnLu} in two dimensions and extend theirs for general non-radial core-shell structures in both two and three dimensions.

\end{abstract}





\tableofcontents

\section{Introduction}
Negative index materials (NIMs) were first investigated theoretically by Veselago in \cite{Veselago} and were innovated by  Nicorovici et al. in \cite{NicoroviciMcPhedranMilton94} and Pendry in~\cite{PendryNegative}. The  existence of such materials was confirmed by Shelby et al.  in \cite{ShelbySmithSchultz}. The study of NIMs has attracted a lot attention in the scientific community thanks to their many applications. One of the appealing ones is cloaking. There are at least three ways to do cloaking using NIMs. The first one is based on plasmonic structures introduced by Alu and Engheta in \cite{AluEngheta}. The second one uses the concept of complementary media. This was suggested by Lai et al. in \cite{LaiChenZhangChanComplementary} and confirmed theoretically in  \cite{Ng-Negative-Cloaking} for a slightly different scheme. 
The last one  is based on the concept of anomalous localized resonance discovered by Milton and Nicorovici in \cite{MiltonNicorovici}. In this paper, we concentrate on the last method.  

\medskip
Cloaking via anomalous localized resonance (CALR) was discovered by Milton and 
 Nicoro-\\vici in \cite{MiltonNicorovici}. Their work has root from \cite{NicoroviciMcPhedranMilton94} (see also \cite{NicoroviciMcPhedranMiltonPodolskiy1}) where the localized resonance was observed and established for constant symmetric plasmonic structures  in the two dimensional quasistatic regime. 
More precisely, in \cite{MiltonNicorovici}, the authors studied  core-shell plasmonic structures  in which  a circular shell has permittivity $-1 + i \delta$ while the core and the matrix, the complement of the core and the shell,  have permittivity 1. Here $\delta$ denotes the loss of the material in the shell.  Let $r_e$ and $r_i$ be the outer and the inner radius of the shell. 
They showed that there is a critical radius  $r_* := (r_e^3 r_i^{-1})^{1/2}$ such that  
a dipole is not seen by an observer away from the core-shell structure, hence it is cloaked,  if and only if the dipole is within distance $r_*$ of the shell; moreover, 
the power $E_\delta(u_\delta)$ of the field $u_\delta$, which is roughly speaking $\delta \| u_\delta\|_{H^{1}}^2$, blows up. 
They called this phenomenon  cloaking via anomalous localized resonance. Two key features of this phenomenon are: 
\begin{enumerate}
\item[1)] the localized resonance, i.e., the fields blow up in some regions and remain bounded in some others as the loss goes to 0. 
\item[2)] the connection between the localized resonance and the blow up of the power as the loss goes to 0. 
\end{enumerate}
Their work has opened a new way of cloaking   and has been a source of inspiration for many investigations see \cite{AmmariCiraoloKangLeeMilton, AmmariCiraoloKangLeeMilton2, AmmariCiraoloKangLeeMilton1, BouchitteSchweizer10, BrunoLintner07,  Chung, KohnLu, KettunenLassas,  Milton-folded, MiltonOnofrei,  MinhLoc, NPET}.


\medskip

Let us discuss recent progress on CALR. In  \cite{BouchitteSchweizer10}, Bouchitte and Schweizer  proved that a small circular inclusion of radius $\gamma(\delta)$ (with $\gamma(\delta) \to 0$ fast enough) is cloaked by the core-shell plasmonic structure  mentioned above in the two dimensional quasistatic regime  if the inclusion is located within distance $r_*$ of the shell. Otherwise it is visible. 
Concerning the second feature of CALR, the blow up of the power  was studied for a more general setting by Ammari et al. in \cite{AmmariCiraoloKangLeeMilton} and  Kohn et al. in \cite{KohnLu}. More precisely, they considered  non-radial core-shell structures in which the shell has permittivity $-1 + i \delta$ and the core and the matrix have permitivity 1. In \cite{AmmariCiraoloKangLeeMilton},  Ammari et al. dealt with   arbitrary shells  in the  two dimensional quasistatic regime. They provided a characterization  of sources for which the power blows up. Their characterization is based on the spectrum of a self-adjoint compact operator (Neumann-Poincar\'e type operator).  In \cite{KohnLu},  Kohn et al. considered core-shell structures  in the two dimensional quasistatic regime in which the matrix is radial symmetric but the core is not.  Using  a variational approach,  they established the blow up of the power 
for a class of sources concentrated on circles within distance $r_* = (r_e^3 r_i^{-1})^{1/2}$ of the core-shell region $ B_{r_e}$ if the core is inside $B_{r_i}$. They also showed that the power remains bounded for a class of sources concentrated on circles outside $B_{r_*}$  if the core is round, inside, and close to $B_{r_i}$. The localized resonance associated with CALR has been so far discussed only for simple geometries, see \cite{AmmariCiraoloKangLeeMilton, AmmariCiraoloKangLeeMilton1, Chung}. 

\medskip

An important class of NIMs in which the localized resonance might appear is the class of reflecting complementary media see \cite{Ng-Negative-Cloaking, Ng-superlensing, MinhLoc1}. The concept of  reflecting complementary media for a general core-shell structure was introduced and studied in \cite{Ng-Complementary}. This class is inspired from the pivotal work of Nicorovici et al. in \cite{NicoroviciMcPhedranMilton94} and from the important  notion of complementary media suggested by  Ramakrishna and Pendry in \cite{PendryRamakrishna0}.  Nevertheless, the complementary property is not enough to ensure that CALR takes place as discussed in \cite{MinhLoc}. Therefore, the study of the two features 1) and 2)  together in CALR is of necessity and importance.

\medskip 
In this paper, we investigate CALR for a subclass of complementary media called 
the class of  doubly complementary media for a core-shell structure,  which will be given in Definition~\ref{def-DCM}.  This class is rich enough to allow us to cloak  an {\bf arbitrary source} concentrating on an {\bf arbitrary smooth bounded manifold of codimension 1} placed in an {\bf  arbitrary medium} via anomalous localized resonance (see Section~\ref{sect-cloaking}); the cloak is independent of the source. Roughly speaking,  the shell is not only  reflecting complementary to a part of the matrix but  also to a part of  the core.   We establish the following three properties on CALR for doubly complementary media, which are what one would expect from a structure for which CALR takes place: 
\begin{enumerate}
\item[P1)] CALR appears if and only if the power blows up (Theorem~\ref{thm1}).
\item[P2)] The power blows up if the source is located  ``near" the shell (Theorem~\ref{thm2}).
\item[P3)] The power remains bounded if the source is far away from the shell (Theorem~\ref{thm3}).
\end{enumerate}
{\bf Property P2)}, the blow up of the power,  is in fact established for {\bf reflecting complementary media}. 
We also address qualitative estimates on the distance from the source to the shell for which CALR does or does not appear  in various situations (Theorems~\ref{thm2} and \ref{thm3}).  




\medskip 
We next describe the problem more precisely.  Let $d=2, 3$, and $\Omega$ be a smooth open bounded subset of $\mR^d$,  and let  $0 < r_1 < r_2$  be such that $B_{r_2} \subset \subset \Omega$. 
Set, for $\delta>0$,  
\begin{equation}\label{def-sd}
s_\delta  : = \left\{\begin{array}{cl} -1 + i \delta & \mbox{ in } B_{r_2} \setminus B_{r_1}, \\[6pt]
1 & \mbox{ otherwise.}
\end{array} \right.
\end{equation}
Let $A$ be a symmetric uniformly elliptic matrix-valued function defined in $\Omega$, i.e., $A$ is symmetric and 
\begin{equation}\label{pro-A}
\frac{1}{\Lambda} |\xi|^2 \le \langle A(x) \xi, \xi \rangle \le \Lambda |\xi|^2 \quad \forall \, \xi \in \mR^d,
\end{equation}
for a.e. $x \in \Omega$ and  for some $1 \le \Lambda < + \infty$. Let $f \in L^2(\Omega)$ with $\supp f \cap B_{r_2} = \O$ and let  
$u_\delta \in H^1_0(\Omega)$ be the unique solution to 
\begin{equation}\label{def-ud}
\dive (s_\delta A \nabla u_\delta) = f \mbox{ in } \Omega. 
\end{equation}
The power $E_\delta(u_\delta)$ is defined by (see, e.g., \cite{MiltonNicorovici})
\begin{equation*}
E_\delta (u_\delta) = \delta  \int_{B_{r_2} \setminus B_{r_1}} |\nabla u_\delta|^2. 
\end{equation*}
Using the fact that $u_\delta = 0$ on $\partial \Omega$, one has \footnote{One way to obtain this inequality is to 
multiply \eqref{def-ud} by $\bar u_\delta$ (the conjugate of $u_\delta$), integrate on $\Omega$, and consider the real part.}
\begin{equation}\label{PPP}
 \int_{\Omega} (|\nabla u_\delta|^2 + |u_\delta|^2) \le C \Big( \int_{B_{r_2} \setminus B_{r_1}} |\nabla u_\delta|^2 + \| f\|_{L^2}^2 \Big), 
\end{equation}
for some positive constants $C$ independent of $f$ and $\delta \in (0, 1)$. Let  $v_\delta \in H^1_0(\Omega)$ be
the unique solution to
\begin{equation}\label{def-vd}
\dive (s_\delta A \nabla v_\delta) = f_\delta \mbox{ in } \Omega. 
\end{equation}
Here 
\begin{equation*}
f_\delta = c_\delta f, 
\end{equation*}
 and  $c_\delta$ is the normalization  constant such that 
\begin{equation}\label{Power-vn}
\delta^{1/2} \int_{B_{r_2} \setminus B_{r_1}} |\nabla v_\delta|^2 = 1. 
\end{equation}

In this paper, we are interested in a class of matrices $A$, the class of doubly complementary media,  for which CALR takes place. 
Before giving the definition of doubly complementary media for a general core-shell structure, let us recall the definition of reflecting complementary media introduced  in \cite[Definition 1]{Ng-Complementary}. 

\begin{definition}[Reflecting complementary media]   \fontfamily{m} \selectfont
 \label{def-Geo} Let $r_1< r_2 < r_3$. The media $A$ in $B_{r_3} \setminus B_{r_2}$ and $-A$ in $B_{r_2} \setminus B_{r_1}$  are said to be  {\it reflecting complementary} if 
there exists a diffeomorphism $F:B_{r_2} \setminus \bar B_{r_1} \to B_{r_3} \setminus \bar B_{r_2}$ such that 
\begin{equation}\label{cond-ASigma}
F_*A  = A  \mbox{ for  } x \in  B_{r_3} \setminus \bar B_{r_2}, 
\end{equation}
\begin{equation}\label{cond-F-boundary}
F(x) = x \mbox{ on } \partial B_{r_2}, 
\end{equation}
and the following two conditions hold:
\begin{enumerate}
\item[1.] There exists an diffeomorphism extension of $F$, which is still denoted by  $F$, from $B_{r_2} \setminus \{x_{1}\} \to \mR^d \setminus \bar B_{r_2}$ for some $x_{1} \in B_{r_1}$. 

\item[2.] There exists a diffeomorphism $G: \mR^d \setminus \bar B_{r_3} \to B_{r_3} \setminus \{x_{2}\}$ for some $x_{2} \in B_{r_3}$  such that 
\begin{equation}\label{cond-G-boundary}
\quad G(x) = x \mbox{ on } \partial B_{r_3},
\end{equation}
and
\begin{equation}\label{extension}
G \circ F : B_{r_1}  \to B_{r_3} \mbox{ is a diffeomorphism if one sets } G\circ F(x_1) = x_2.
\end{equation}
\end{enumerate} 
\end{definition}

Here and in what follows, if $T$ is a diffeomorphism and  $a$ is a matrix-valued function, we denote  
\begin{equation}\label{def-F*}
T_*a(y) = \frac{DT(x) a(x) DT(x)^T}{|\det DT(x)|} \quad  \mbox{ where } x = T^{-1}(y). 
\end{equation}

\begin{remark}  \fontfamily{m} \selectfont
In \eqref{cond-F-boundary} and \eqref{cond-G-boundary}, $F$ and $G$ denote some  diffeomorphism extensions of $F$ and $G$ in a neighborhood of $\partial B_{r_2}$ and of $\partial B_{r_3}$. As noted in \cite{Ng-Complementary}, conditions  \eqref{cond-ASigma} and \eqref{cond-F-boundary} are the main assumptions in Definition~\ref{def-Geo}. The term ``reflecting'' in Definition~\ref{def-Geo} comes from \eqref{cond-F-boundary} and the fact that $B_{r_1} \subset B_{r_2} \subset B_{r_3}$. Conditions 1) and 2)  are  mild assumptions. Introducing $G$ makes the analysis more accessible, see \cite{Ng-Complementary, Ng-Negative-Cloaking, Ng-superlensing, MinhLoc1} and the analysis presented in this paper. 
\end{remark}

\begin{remark} \fontfamily{m} \selectfont The class of reflecting complementary media has played an important role in the other applications of NIMs such as cloaking and superlensing using complementary see \cite{Ng-Negative-Cloaking, Ng-superlensing, MinhLoc1}. 
\end{remark}

\begin{remark} \fontfamily{m} \selectfont  Taking $d=2$,  $A = I$ and $r_3 = r_2^2/r_1$, and letting $F$ be the Kelvin transform with respect to $\partial B_{r_2}$, i.e., $F(x) = r_2^2 x/ |x|^2$,  one can  verify that the core-shell structures considered by Milton et al. in \cite{MiltonNicorovici} and by Kohn et al. in \cite{KohnLu}  have the reflecting complementary property.  

\end{remark}

We are ready to introduce the concept of  doubly complementary media. 

\begin{figure}[h!]
\begin{center}
\includegraphics[width=8cm]{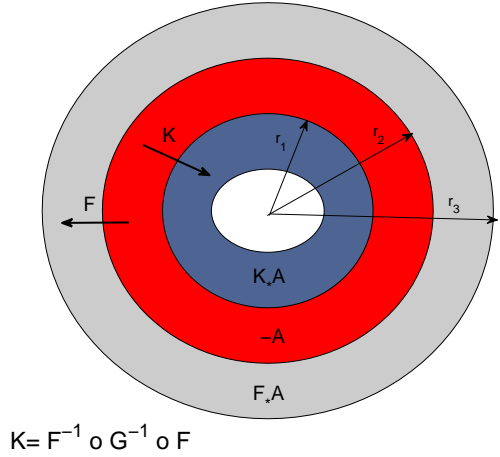}
\caption{$s_0A$ is doubly complementary: $-A$ in $B_{r_2} \setminus B_{r_1}$ (the red region) is complementary to $A = F_*A$ in $B_{r_3} \setminus B_{r_2}$ (the grey region) and $A = K_*A$ with $K = F^{-1} \circ G^{-1} \circ F$ in $K(B_{r_2} \setminus B_{r_1})$ (the blue grey region).} \label{fig1}
\end{center}
\end{figure}

\begin{definition} \label{def-DCM} \fontfamily{m} \selectfont The medium $s_0 A$ is said to be {\it doubly complementary} if  for some $r_3 >0$ with $B_{r_3} \subset \subset \Omega$, $A$ in $B_{r_3} \setminus B_{r_2}$ and $-A$ in $B_{r_2} \setminus B_{r_1}$ are reflecting complementary, and 
\begin{equation}\label{def-DC}
F_*A = G_* F_*A = A \mbox{ in } B_{r_3} \setminus B_{r_2}, 
\end{equation}
for some $F$ and $G$ coming from Definition~\ref{def-Geo} (see Figure \ref{fig1}). 
\end{definition}

\begin{remark}\fontfamily{m} \selectfont
The reason for which media satisfying \eqref{def-DC} are called doubly complementary media is that $-A$ in $B_{r_2} \setminus B_{r_1}$ is not only complementary to $A$ in $B_{r_3} \setminus B_{r_2}$ but also to $A$ in $(G \circ F)^{-1}(B_{r_3} \setminus \overline{B_{r_2}})$ (a subset of $B_{r_1}$) (see \cite{Ng-CALR-CRAS}). 
\end{remark}

\begin{remark}\fontfamily{m} \selectfont
Taking $d=2$, $A = I$ and $r_3 = r_2^2/r_1$, and letting $F$ and $G$ be the Kelvin transform with respect to $\partial B_{r_2}$ and $\partial B_{r_3}$,  one can  verify that the core-shell structures considered by Milton et al. in \cite{MiltonNicorovici} is of doubly complementary property.  The setting considered in \cite{AmmariCiraoloKangLeeMilton2} also has this property. 
\end{remark}

In what follows, we assume that 
\begin{equation}\label{AC3}
A \in [C^3(\overline{B_{r_3} \setminus B_{r_2}})]^{d \times d}. 
\end{equation}
This assumption, which  can be sometimes weaken,  is necessary for the use of a three sphere inequality, the unique continuation principle, and  the separation of variables technique introduced later in this paper.

\medskip
The following theorem is one of the main results of the paper. It is on the equivalence between the blow up of the power and CALR for doubly complementary media,  which  implies Property P1).

\begin{theorem}\label{thm1} Let $d=2,  3$, $f \in L^2(\Omega)$ with $\supp f \subset \Omega \setminus B_{r_2}$, $(\delta_n) \to0$, and let $u_{\delta_n} \in H^1_0(\Omega)$ be the unique solution to 
\begin{equation*}
\dive (s_{\delta_n} A \nabla u_{\delta_n}) = f \mbox{ in } \Omega. 
\end{equation*}
Assume that $s_0A$ is doubly complementary. 
We have 
\begin{itemize}
\item[i)] If $\lim_{n \to \infty} \delta_n \| \nabla u_{\delta_n}\|_{L^2(B_{r_2} \setminus B_{r_1})}^2 = + \infty$, then 
\begin{equation}\label{LCR-vn}
v_{\delta_n} \to 0 \mbox{ weakly in } H^1(\Omega \setminus B_{r_3}), 
\end{equation}
where $v_\delta \in H^1_0(\Omega)$ is defined in \eqref{def-vd}. 

\item[ii)] If $ \big(\delta_n\| \nabla u_{\delta_n}\|_{L^2(B_{r_2} \setminus B_{r_1})}^2 \big)_{n \in \mN}$  is  bounded then 
\begin{equation*}
u_{\delta_n} \to u \mbox{ weakly in } H^1(\Omega \setminus B_{r_3}) \mbox{ as } \delta \to 0, 
\end{equation*}
where $u \in H^1_0(\Omega)$ is the unique solution to 
\begin{equation}\label{UU}
\dive(\hat A \nabla u) = f \mbox{ in } \Omega. 
\end{equation}

\end{itemize}
\end{theorem}

Here and in what follows, we denote
\begin{equation}\label{def-hatA}
\hat A = \left\{\begin{array}{cl} A & \mbox{ in } \Omega \setminus B_{r_3}, \\[6pt]
G_*F_*A   & \mbox{ in } B_{r_3}. 
\end{array} \right.
\end{equation}

The proof  of Theorem~\ref{thm1} is given in Section~\ref{sect-thm1} where a stronger result (Proposition~\ref{pro1.1}) is established. 

\medskip
The equivalence between the blow up of the power and the CALR can be obtained from  Theorem~\ref{thm1} as follows. Suppose that the power blows up,  i.e.,  
\begin{equation*}
\lim_{n \to \infty} \delta_n \| \nabla u_{\delta_n}\|_{L^2(B_{r_2} \setminus B_{r_1})}^2 = + \infty. 
\end{equation*}
Then, by Theorem~\ref{thm1},  $v_{\delta_n} \to 0$ in $\Omega \setminus B_{r_3}$. The source $\alpha_{\delta_n} f$ is not seen by observers far away from the shell: the source is cloaked.  We note that the localized resonance happens   
in this case since both \eqref{Power-vn}  and \eqref{LCR-vn}  take place. If the power of $u_{\delta_n}$ remains bounded,  then $u_{\delta_n} \to u$ weakly in $H^1(\Omega \setminus B_{r_3})$. Since $u \in H^1_0(\Omega)$ is the unique solution to \eqref{UU}, the source is not cloaked.  

\medskip 
Theorem~\ref{thm1} is, to our knowledge, the first result providing the connection between the blow up of the power and the invisibility of a source in  a  general setting. The  standard separation of variables is out of reach here.

\medskip
We next show that CALR takes place if the source is  located ``near" the shell. This implies Property P2). In fact, we establish this property for  {\bf reflecting  complementary media}.  More precisely, we have the following result whose proof is given in Section~\ref{sect-thm2}.  

\begin{theorem}\label{thm2} Let $d=2, 3$, $f \in L^2(\Omega)$ with $\supp f \subset \Omega \setminus B_{r_2}$, and let $u_\delta \in H^1_0(\Omega)$ be the unique solution to 
\begin{equation*}
\dive(s_\delta A \nabla u_\delta) = f \mbox{ in } \Omega. 
\end{equation*}
Assume that $A$ in $B_{\hat r_3} \setminus B_{r_2}$ and $-A$ in $B_{r_2} \setminus B_{\hat r_1}$ are reflecting complementary for some $r_1 \le \hat r_1 <  r_2 <  \hat r_3$,  with $B_{\hat r_3} \subset \subset \Omega$. 
There exists a constant $r_* \in (r_2, \hat r_3)$, independent of $\delta$ and $f$ such that if there is {\bf no}  $w \in H^1(B_{r_*} \setminus B_{r_2})$  with the properties 
\begin{equation}\label{pro-BU}
\dive (A \nabla w) = f \mbox{ in } B_{r_*} \setminus B_{r_2}, \quad w = 0 \mbox{ on } \partial B_{r_2}, \quad \mbox{ and } \quad A \nabla w \cdot \eta = 0  \mbox{ on } \partial B_{r_2}, 
\end{equation}
then 
\begin{equation}\label{blow-up}
\limsup_{\delta \to 0} \delta^{1/2} \| \nabla u_\delta \|_{L^2(B_{r_2} \setminus B_{r_1})} = + \infty. 
\end{equation}
Assume in addition that $A = I$ in $B_{\hat r_3} \setminus B_{r_2}$, then 
\begin{equation}\label{pro-BU-0}
r_* \mbox{ can be taken by any number less than } \sqrt{\hat r_3 r_2}. 
\end{equation}
\end{theorem}

Here and in what follows, for $D$ a smooth bounded open subset of $\mR^d$, on $\partial D$,  $\eta$ denotes the outward unit normal vector.

\medskip
Concerning the boundedness of the power, we prove
\begin{theorem}\label{thm3}
Let $d=2, \, 3$, $f \in L^2(\Omega)$, and let $u_\delta \in H_0^1(\Omega)$ be the unique solution to \eqref{def-ud}. Assume that $s_0A$  is doubly complementary and $\supp f \cap B_{r_3} = \O$. Then 
\begin{equation}\label{part1-thm3}
\limsup_{\delta \to 0} \| u_\delta\|_{H^1(\Omega)} < + \infty. 
\end{equation}
Assume in addition that $A = I$ in $B_{r_3} \setminus B_{r_2}$. 
If  there exists $w \in H^1(B_{r_0} \setminus B_{r_2})$ for some $r_0 > \sqrt{r_2 r_3}$ with the properties 
\begin{equation*}
\dive (A \nabla w) = f \mbox{ in } B_{r_0} \setminus B_{r_2}, \quad w = 0 \mbox{ on } \partial B_{r_2}, \quad \mbox{ and } \quad  A \nabla w \cdot \eta = 0 \mbox{ on } \partial B_{r_2}, 
\end{equation*}
then 
\begin{equation}\label{part2-thm3}
\limsup_{\delta \to 0} \delta^{1/2} \|u_\delta \|_{H^{1}(\Omega)} < + \infty. 
\end{equation}
\end{theorem}

It is clear that Theorem~\ref{thm3} implies Property P3). The proof of Theorem~\ref{thm3} is given in Section~\ref{sect-thm3}. 

\medskip

The analysis in this paper is based on several new observations and ideas. The proof of Theorem~\ref{thm1} (in Section~\ref{sect-thm1}) make uses of  the reflecting and the removing localized singularity techniques introduced in \cite{Ng-Complementary, Ng-Negative-Cloaking, Ng-superlensing} to deal with the localized resonance. To develop these techniques 
for a general core-shell structure, we introduce and implement the separation of variables technique to solve Cauchy problems in a general shell (Proposition~\ref{pro1} in Section~\ref{sect-separation}). 
 The way to implement this technique is one of the cores of the analysis in this paper. 
The use of separation of variables to solve boundary value problems for the Laplace equation in an arbitrary domain was considered in the literature and was based on  the integral method, see e.g.,  \cite{KarpShamma}. The analysis presented here is based  on the idea of transformation optics and the reflecting technique.  As a consequence, we obtain the existence of surface plasmons for general complementary media (Proposition~\ref{pro1}). 
The proof of Theorem~\ref{thm2} (in Section~\ref{sect-thm2}) is based on a new observation for complementary media (Lemma~\ref{lem-pro2.1}) whose proof is based on a three spheres inequality. The idea of the proof of  Theorem~\ref{thm3} (in Section~\ref{sect-thm3})  is as follows. The first part \eqref{part1-thm3} is from \cite{Ng-Complementary}. The proof of the second part \eqref{part2-thm3}  is based on a  kind of  removing singularity technique and  uses ideas in \cite{Ng-Complementary}. A key point is the construction of an auxiliary function  $W_\delta$ in \eqref{def-Vd-S}.  Using Theorem~\ref{thm1} and Theorem~\ref{thm2}, we can construct a cloaking device to cloak a general source concentrate on a manifold of codimension 1 in an arbitrary medium (see Section~\ref{sect-cloaking}). The proof also makes use of  the unique continuation principle. 

 \medskip By considering $A=I$ in Theorems~\ref{thm1}, \ref{thm2}, and  \ref{thm3}, one can rediscover the results of Milton and Nicorovici in  \cite{MiltonNicorovici},  and Kohn et al. in \cite{KohnLu},  and the results of Ammari et al.  in \cite{AmmariCiraoloKangLeeMilton}  in the radial setting,  in two dimensions.  The results presented here extend theirs for general non-radial core-shell structures in both two and three dimensions. 

\medskip

The results in this paper are announced in \cite{Ng-CALR-CRAS}. The study of CALR  in the finite frequency regime will be considered in \cite{Ng-CALR-finite}.

\section{A condition on the blow up of the power. Proof of Theorem~\ref{thm2}}\label{sect-thm2}

This section containing two subsections is devoted to the proof of Theorem~\ref{thm2}. 
 In the first subsection, we present some useful lemmas.  
 The proof of Theorem~\ref{thm2} is given in  the second subsection. 
 
\subsection{Preliminaries}\label{sect-pre-thm2}

We first recall the following result, a change of variables formula, which follows immediately from 
\cite[Lemma 2]{Ng-Complementary}, and 
is used repeatedly in this paper.

\begin{lemma}\label{lem-TO} Let $d=2,3$, $R>0$, $D_1$ and $D_2$ be two smooth open subsets of $\mR^d$ be such that $D_1 \subset \subset B_R  \subset \subset D_2$. Assume that  $T$ is a diffeomorphism from $B_R \setminus D_1$ onto $D_2 \setminus B_R$ and let $a \in [L^\infty(B_R \setminus D_1)]^{d \times d}$ be uniformly elliptic. Fix $u \in H^1(B_R \setminus D_1)$ and set $v = u \circ T^{-1}$. Then
\begin{equation*}
\dive (a \nabla u)  = 0 \mbox{ in } B_R \setminus D_{1}
\end{equation*}
if and only if
\begin{equation*}
\dive (T_*a \nabla v) = 0\mbox{ in } D_{2} \setminus B_R. 
\end{equation*}
Assume in addition that $T(x) = x$ on $\partial B_R$.  We have
\begin{equation}\label{reflextion}
T_*a \nabla v \cdot \eta = - a \nabla u \cdot \eta  \mbox{ on } \partial B_R.
\end{equation}
\end{lemma}

We next recall the following known result on three spheres inequalities (see, e.g.,  \cite[Theorem 2.3 and (2.10)]{AlessandriniRondi}). 

\begin{lemma} [Three spheres inequality] \label{lem1.1}  Let $d =2, 3$, $0< R_{1} < R_{2} < R_{3}$, and let $M$ be a Lipschitz matrix-valued function defined in $B_{R_{3}}$ such that  $M$ is symmetric and  uniformly elliptic in $B_{R_3}$,  and $M(0) = I$. 
Assume  $v \in H^{1}(B_{R_{3}})$ is a solution to 
\begin{equation*}
\dive (M \nabla v ) = 0 \mbox{ in } B_{R_{3}}. 
\end{equation*}
There exist two positive constants $C$ and $c$   depending only on $R_3$ and the ellipticity  and the Lipschitz constants of $M$ such that 
\begin{equation*}
\|v \|_{L^2(\partial B_{R_{2}} )} \le C \| v\|_{L^2(B_{\partial R_{1}})}^{ \alpha }  \| v\|_{L^2(\partial B_{R_{3}})}^{1 - \alpha }, 
\end{equation*}
where 
\begin{equation}\label{alpha-1}
\alpha = \frac{\ln ( \frac{R_{3}}{R_{2}}) } { \ln (\frac{R_{3}}{R_{2}}) + c \ln (\frac{R_{2}}{ R_{1}})}, 
\end{equation}
In the case $M = I$ in $B_{R_3}$,  one can take $c=1$,  i.e., 
\begin{equation*}
\alpha = \ln \Big(\frac{R_3}{R_2} \Big) \Big/ \ln \Big( \frac{R_3}{R_1} \Big). 
\end{equation*}
\end{lemma}

Using  Lemma~\ref{lem1.1}, we can prove

\begin{lemma} \label{lem1.2}  Let $d = 2, 3$,  $0< R_{1} < R_{2} < R_{3}$,  and let $M$ be a Lipschitz matrix-valued function defined in $B_{R_{3}}$ such that  $M$ is symmetric and uniformly elliptic in $B_{R_3}$ and $M(0) = I$. 
Assume  $v \in H^{1}(B_{R_{3}})$ is a solution to 
\begin{equation*}
\dive (M \nabla v ) = 0 \mbox{ in } B_{R_{3}} \setminus B_{R_1}. 
\end{equation*}
There exist two positive constants $C$ and  $c$ such that $C$ depends only on $R_1, R_3$, the ellipticity  and  the Lipschitz constants of $M$, and $c$ depends only on $R_3$, the ellipticity  and  the Lipschitz constants of $M$,  and 
\begin{multline}\label{SI1}
\|v \|_{L^2(\partial B_{R_2})} \le C \left( \big(\| v\|_{H^{1/2}(B_{\partial R_{1}})}+ \| M \nabla  v \cdot \eta\|_{H^{-1/2}(B_{\partial R_{1}})}  \big)^{ \alpha }  \| v\|_{L^2(\partial B_{R_{3}})}^{1 - \alpha }  \right.\\[6pt]+ 
\left. \big(\| v\|_{H^{1/2}(B_{\partial R_{1}})}+ \| M \nabla  v \cdot \eta\|_{H^{-1/2}(B_{\partial R_{1}})}  \big) \right), 
\end{multline}
where 
\begin{equation}\label{alpha-2}
\alpha = \frac{\ln ( \frac{R_{3}}{R_{2}}) } { \ln (\frac{R_{3}}{ R_{2}}) + c \ln (\frac{ R_{2}}{R_{1}})}. 
\end{equation}
In the case $M = I$ in $B_{R_3}$,  one can take $c=1$. 
\end{lemma}

\noindent{\bf Proof.} Let $w \in H^1(B_{R_3} \setminus \partial B_{R_1})$ be such that 
\begin{equation*}
\dive (M \nabla w) = 0 \mbox{ in } B_{R_3} \setminus \partial B_{R_1}, \quad w = 0 \mbox{ on } \partial B_{R_3}, 
\end{equation*}
\begin{equation*}
[w]= v   \mbox{ on } \partial B_{R_1},  \quad \mbox{ and }  \quad [M \nabla w \cdot \eta]  = M \nabla v \cdot \eta \mbox{ on } \partial B_{R_1}.  
\end{equation*}
Henceforth   $[\cdot ]$ denotes the jump across the boundary. 
It follows that 
\begin{equation}\label{tototototo}
\|w \|_{H^1(B_{R_3} \setminus \partial B_{R_1})} \le C \big(  \| v \|_{H^{1/2}(\partial B_{R_1})} + \| M \nabla v \cdot \eta  \|_{H^{-1/2}(\partial B_{R_1})}\big). 
\end{equation}
Here and in what follows in this proof, $C$ denotes a positive constant depending only on $R_1$, $R_3$, and the ellipticity  and  the Lipschitz constants of $M$.  Define 
\begin{equation*}
V = \left\{\begin{array}{cl} v - w & \mbox{ in } B_{R_3} \setminus B_{R_1}, \\[6pt]
- w & \mbox{ in } B_{R_1}. 
\end{array}\right.
\end{equation*}
Then $V \in H^1(B_{R_3})$ and 
\begin{equation*}
\dive (M \nabla V) = 0 \mbox{ in } B_{R_3}. 
\end{equation*}
Applying Lemma~\ref{lem1.1}, we have
\begin{equation*}
\| V \|_{L^2(\partial B_{R_2})} \le C \| V \|_{L^2(\partial B_{R_1})}^\alpha \| V \|_{L^2(\partial B_{R_3})}^{1 - \alpha}. 
\end{equation*}
The conclusion follows from \eqref{tototototo} and the definition of $V$. \proofend

\medskip
The following result provides the key ingredient for the proof of  Theorem~\ref{thm2}.

\begin{lemma}\label{lem-pro2.1} Let $d=2, 3$,  $0 < R_1 < R_2< + \infty$, $M$ be a symmetric uniformly elliptic matrix-valued function defined  in $B_{R_2} \setminus B_{R_1}$, and let 
$g, h \in L^2(B_{R_2} \setminus B_{R_1})$. Assume that $M$ is Lipschitz  and $U_\delta, V_\delta  \in H^1(B_{R_2} \setminus B_{R_1})$ satisfy
\begin{equation*}
\dive(M \nabla U_\delta) = g \mbox{ in } B_{R_2} \setminus B_{R_1}, \quad  \dive(M \nabla V_\delta) = h \mbox{ in } B_{R_2} \setminus B_{R_1}, 
\end{equation*}
\begin{equation*}
U_\delta  = V_\delta \mbox{ on } \partial B_{R_1}, \quad \mbox{ and } \quad  M \nabla U_\delta \cdot \eta  =  (1 - i \delta) M \nabla V_\delta \cdot \eta \mbox{ on } \partial  B_{R_1}. 
\end{equation*}
There exists a constant $R_* \in (R_1, R_2)$ depending only on $R_1$, $R_2$, and the ellipticity and the Lipschitz  constants of $M$, but  independent of $\delta$, $g$, and $h$ such that if there is {\bf no}  $W \in H^1(B_{R_*} \setminus B_{R_1})$ with the properties 
\begin{equation}\label{pro-BU-1}
\dive (M \nabla W) = g- h \mbox{ in } B_{R_*} \setminus B_{R_1}, \quad W = 0 \mbox{ on } \partial B_{R_1}, \quad \mbox{and} \quad M \nabla W \cdot \eta = 0  \mbox{ on } \partial B_{R_1}, 
\end{equation}
then 
\begin{equation}\label{blow-up-1}
\limsup_{\delta \to 0} \delta^{1/2} \Big( \| U_\delta\|_{H^1(B_{R_2} \setminus B_{R_1})}   + \| V_\delta \|_{H^1(B_{R_2} \setminus B_{R_1})} \Big)= + \infty. 
\end{equation}
Assume in addition that $M = I$ in $B_{R_2} \setminus B_{R_1}$, then 
\begin{equation}\label{pro-BU-1-1}
R_* \mbox{ can be taken by any number less than } \sqrt{R_1 R_2}. 
\end{equation}
\end{lemma}

\noindent{\bf Proof.} For notational ease, we denote $U_{2^{-n}}$ and $V_{2^{-n}}$ by  $U_n$ and $V_n$. We have
\begin{equation*}
\dive(M \nabla U_n) = g \mbox{ in } B_{R_2} \setminus B_{R_1}, \quad  \dive(M \nabla V_n) = h \mbox{ in } B_{R_2} \setminus B_{R_1}, 
\end{equation*}
and 
\begin{equation*}
U_n  = V_n \mbox{ on } \partial B_{R_1}, \quad  M \nabla U_n \cdot \eta  =  (1 - i 2^{-n}) M \nabla V_n \cdot \eta \mbox{ on } \partial  B_{R_1}.
\end{equation*}
Let  $\hat M$ be an extension on $M$ in $B_{R_2}$ such that $\hat M$ is Lipschitz and uniformly elliptic in $B_{R_2}$, and $\hat M(0) = I$ \footnote{One can choose $\hat M$ as follows:  $\hat M(x) = (2 r/ R_1 -1 ) M(R_1 \sigma) + (2 - 2 r/ R_1) I $ if $x \in B_{R_1} \setminus B_{R_1/2}$ and $\hat M(x)  = I$ if $x \in B_{R_1/2}$, where $r = |x|$ and $\sigma = x/ |x|$.  In the case $M=I$ in $B_{R_2} \setminus B_{R_1}$, we choose $\hat M = I$ in $B_{R_1}$. }.  Let $c$ be the constant in Lemma~\ref{lem1.2} corresponding to $\hat M$ and the shell $B_{R_2} \setminus B_{R_1}$. Define
\begin{equation*}
\alpha(r) =  \frac{\ln ( \frac{R_{2}}{r}) } { \ln (\frac{R_{2}}{r}) + c \ln (\frac{r}{R_{1}})} \quad \forall \, r \in (R_1, R_2). 
\end{equation*}
Fix $R_*$ such that $\alpha(R_*) > 1/2$ (this holds if $R_*$ is chosen close to $R_1$). There exists $\gamma \in (0, 1)$ (close to 1) such that 
\begin{equation}\label{choice-r}
\alpha(r) > \big(\alpha(R_*) + 1/2 \big)/2 \quad \mbox{ for } r \in (\gamma R_*, (2 - \gamma) R_*). 
\end{equation}

We prove by contradiction that
\begin{equation}\label{claim-IE}
\limsup_{n \to +\infty} 2^{-n/2} \big( \| U_n\|_{H^{1}(B_{R_2} \setminus B_{R_1})}   + \| V_n \|_{H^{1}(B_{R_2} \setminus B_{R_1})}  \big)= + \infty. 
\end{equation} 
Assume that 
\begin{equation}\label{claim-IE-contradiction}
m: = \sup_{n}  2^{-n/2} \big( \| U_n\|_{H^{1}(B_{R_2} \setminus B_{R_1})}   + \| V_n \|_{H^{1}(B_{R_2} \setminus B_{R_1})}  \big)  <  + \infty. 
\end{equation} 
Define 
\begin{equation*}
W_n = U_n- V_n \mbox{ in } B_{R_2} \setminus B_{R_1} \quad \mbox{ and } \quad \Phi_n = -  i 2^{-n} M \nabla V_n \cdot \eta \mbox{ on } \partial B_{R_1}. 
\end{equation*}
Then 
\begin{equation*}
\dive(M \nabla W_n) = g - h \mbox{ in } B_{R_2} \setminus B_{R_1}, \; \;  W_n = 0 \mbox{ on } \partial B_{R_1}, \mbox{  and  } \quad M \nabla W_n \cdot \eta =  \Phi_n \mbox{ on } \partial B_{R_1}. 
\end{equation*}

We claim that $(W_n)$ is a Cauchy sequence in $H^1(B_{R_*} \setminus B_{R_1})$. 

Indeed, set 
\begin{equation*}
w_{n} = W_{n+1} - W_{n}  \mbox{ in } B_{R_2} \setminus B_{R_1} \quad \mbox{ and } \quad \phi_{n} = \Phi_{n+ 1} - \Phi_{n} \mbox{ on } \partial B_{R_1}. 
\end{equation*}
We have
\begin{equation*}
\dive(M \nabla w_{n}) = 0 \mbox{ in } B_{R_2} \setminus B_{R_1}, \quad w_{n}= 0 \mbox{ on } \partial B_{R_1}, \quad \mbox{ and } \quad A \nabla w_{n} \cdot \eta = \phi_{n} \mbox{ on } \partial B_{R_1}. 
\end{equation*}
From \eqref{claim-IE-contradiction}, we derive that 
\begin{equation*}
\| w_n\|_{H^1(B_{R_2} \setminus B_{R_1})} \le C m 2^{n/2} \quad \mbox{ and } \quad \| \phi_n\|_{H^{1/2}(\partial B_{R_1})} \le C m 2^{-n/2}. 
\end{equation*}
In this proof, $C$ denotes a constant independent of $n$. 
Applying Lemma~\ref{lem1.2}, we obtain 
\begin{equation*}
\| w_{n} \|_{L^2(\partial B_{r})} \le C \Big( \| \phi_n \|_{H^{-1/2}(\partial B_{R_1}) }^{\alpha(r)} \| w_{n}\|_{L^2(\partial B_{R_2})}^{1 - \alpha(r)} + \| \phi_n \|_{H^{-1/2}(\partial B_{R_1}) } \Big) \le C m 2^{-n \beta(r)},  
\end{equation*}
where 
\begin{equation*}
\beta(r) = \big(2 \alpha(r) - 1 \big)/ 2.  
\end{equation*}
From  \eqref{choice-r},
\begin{equation*}
\beta(r)>  \big(\alpha(R_*) - 1/2 \big)/2 > 0 \quad \mbox{ for } r \in (\gamma R_*, (2 - \gamma) R_*).
\end{equation*}
Since $\dive (M \nabla w_n) = 0$ in $B_{R_2} \setminus B_{R_1}$, by the regularity theory of elliptic equations, we obtain 
\begin{equation*}
\| w_{n} \|_{H^{1/2}(\partial B_{R_*})} \le  C m 2^{-n \big(\alpha(R_*) - 1/2 \big)/2}. 
\end{equation*}
Since  $\dive (M \nabla w_n) = 0$ in $B_{R_*} \setminus B_{R_1}$ and $w_n = 0$ on $\partial B_{R_1}$, it follows that 
\begin{equation*}
\| w_{n} \|_{H^{1}( B_{R_*} \setminus B_{R_1} )} \le  C m 2^{-n \big(\alpha(R_*) - 1/2\big)/2}. 
\end{equation*}
Hence $(W_n)$ is a Cauchy sequence  in $H^1(B_{R_*} \setminus B_{R_1})$. 
Let $W $ be the  limit of $W_n$ in $H^1(B_{R_*} \setminus B_{R_1})$.  Then  
\begin{equation*}
\dive (M \nabla W) = g- h \mbox{ in } B_{R_*} \setminus B_{R_1}, \quad  W = 0 \mbox{ on } \partial B_{R_1}, \quad  M \nabla W \cdot \eta =  0 \mbox{ on } \partial B_{R_1}. 
\end{equation*}
This contradicts the non-existence of $W$.  Hence \eqref{claim-IE} holds. The proof is complete. \proofend
 
\subsection{Proof of Theorem~\ref{thm2}}

Set 
\begin{equation*}
u_{1, \delta} = u_{\delta} \circ F^{-1} \mbox{ in } B_{\hat r_3} \setminus B_{r_2}.
\end{equation*}
Since $F_*A = A$ in $B_{\hat r_3} \setminus B_{r_2}$ and $F(x) = x$ on $\partial B_{r_2}$, it follows from Lemma~\ref{lem-TO} that 
\begin{equation*}
\dive(A \nabla u_{1, \delta}) = 0 \mbox{ in } B_{\hat r_3} \setminus B_{r_2}, 
\end{equation*}
\begin{equation*}
u_\delta  = u_{1, \delta} \mbox{ on } \partial B_{r_2}, \quad \mbox{ and } \quad   A \nabla u_\delta \cdot \eta  =  (1 - i \delta) A \nabla u_{1, \delta} \cdot \eta \mbox{ on } \partial  B_{r_2}. 
\end{equation*}
Recall that 
\begin{equation*}
\dive(A \nabla u_{\delta}) = f \mbox{ in } B_{\hat r_3} \setminus B_{r_2}. 
\end{equation*}
Applying Lemma~\ref{lem-pro2.1} with $U_\delta = u_\delta$, $V_\delta = u_{1, \delta}$, $R_1= r_2$, and $R_2 = \hat r_3$, there exists  a constant $r_* \in (r_2, r_3)$, independent of $\delta$ and $f$ such that if there is no solution 
$w \in H^1(B_{r_*} \setminus B_{r_2})$  to \eqref{pro-BU},  then 
\begin{equation*}
\limsup_{\delta \to 0} \delta^{1/2} \big( \|u_\delta \|_{H^1(B_{\hat r_3} \setminus B_{r_2})} + \|u_{1, \delta} \|_{H^1(B_{\hat r_3} \setminus B_{r_2})} \big) = + \infty. 
\end{equation*}
This implies, by \eqref{PPP},  
\begin{equation*}
\limsup_{\delta \to 0} \delta^{1/2} \|\nabla u_\delta\|_{L^2(B_{r_2} \setminus B_{r_1})} = + \infty. 
\end{equation*}
In the case $A = I$ in $B_{\hat r_3} \setminus B_{r_2}$, by Lemma~\ref{lem-pro2.1},  $r_*$ can be taken by any number less than $\sqrt{\hat r_3 r_2}$. \proofend

\section{A condition on the boundedness of the power. Proof of Theorem~\ref{thm3}} \label{sect-thm3}

This section contains two subsections. In the first subsection, we present two lemmas used in the proof of Theorem~\ref{thm3}. The proof of Theorem~\ref{thm3} is given in the second subsection.

\subsection{Two useful lemmas}

The first lemma was established in \cite[Lemma 1]{Ng-Complementary}. 

\begin{lemma} \label{lem-stability-Ng} Let $d=2, 3$, $\delta \in (0, 1)$, and $f \in H^{-1}(\Omega)$ and  let $u_\delta \in H^{1}_0(\Omega)$ be the unique solution to 
\begin{equation*}
\dive(s_\delta A \nabla u_\delta) = f \mbox{ in } \Omega. 
\end{equation*}
Then 
\begin{equation*}
\|u_\delta \|_{H^1(\Omega)} \le \frac{C}{\delta} \| f\|_{H^{-1}(\Omega)}, 
\end{equation*}
for some positive constant $C$ independent of $f$ and $\delta$. 
\end{lemma}

Here is the second lemma whose  proof has root from \cite{Ng-Complementary}. 

\begin{lemma}\label{lem-compatible} Let $d=2, 3$, $\delta \in (0, 1)$, let $f \in L^2(\Omega)$, $g \in H^{1/2}(\partial B_{r_3})$,  and $h \in H^{-1/2}(\partial B_{r_3})$. Assume that  $s_0 A$ is doubly complementary and $\supp f \subset \Omega \setminus B_{r_3}$,  and  let $V_\delta \in H^1(\Omega \setminus \partial B_{r_3})$ be the unique solution to 
\begin{equation*}
\left\{\begin{array}{c}\dive (s_\delta A \nabla V_\delta) = f \mbox{ in } \Omega \setminus \partial B_{r_3},  \\[6pt]
[V_\delta] = g \quad \mbox{ and } \quad [A \nabla V_\delta \cdot \eta] = h \mbox{ on } \partial B_{r_3}, \\[6pt]
V_\delta = 0 \mbox{ on } \partial \Omega. 
\end{array}\right.
\end{equation*}
Then 
\begin{equation*}
 \|V_\delta \|_{H^{1}(\Omega \setminus \partial B_{r_3})} \le C \big(\| f\|_{L^2(\Omega)} + \| g\|_{H^{1/2}(\partial B_{r_3})} + \| h\|_{H^{-1/2}(\partial B_{r_3})}\big),  
\end{equation*}
for some positive constant $C$ independent of $\delta$, $f$, $g$, and $h$. 
\end{lemma} 

\begin{remark} \fontfamily{m} \selectfont The case $g= h = 0$ was considered in \cite{Ng-Complementary}  (see \cite[Theorem 1 and Corollary 1]{Ng-Complementary}). 

\end{remark}

\noindent{\bf Proof.} Let $U \in H^1(\Omega \setminus \partial B_{r_3})$ be the unique solution to 
\begin{equation*}
\left\{\begin{array}{c}\dive (\hat A \nabla U) = f \mbox{ in } \Omega \setminus \partial B_{r_3},  \\[6pt]
[U] = g \quad \mbox{ and } \quad [\hat A \nabla U \cdot \eta] = h \mbox{ on } \partial B_{r_3}, \\[6pt]
U = 0 \mbox{ on } \partial \Omega,  
\end{array}\right.
\end{equation*}
where $\hat A$ is defined in \eqref{def-hatA}. 
Then
\begin{equation}\label{haha1-Thm2.3}
\| U \|_{H^1(\Omega \setminus \partial B_{r_3})} \le C \big(\| f\|_{L^2(\Omega)} + \| g\|_{H^{1/2}(\partial B_{r_3}) } + \| h\|_{H^{-1/2}(\partial B_{r_3})}\big).  
\end{equation}
Define $V_0 \in H^1(\Omega \setminus \partial B_{r_3})$ as follows 
\begin{equation}\label{haha2-Thm2.3}
V_0 =  \left\{\begin{array}{cl} U &\mbox{ in } \Omega \setminus B_{r_2},  \\[6pt]
U \circ F & \mbox{ in } B_{r_2} \setminus B_{r_1}, \\[6pt]
U \circ G \circ F & \mbox{ in } B_{r_1}. 
\end{array}\right.
\end{equation}
Using \eqref{def-DC} and applying Lemma~\ref{lem-TO}, as in \cite[Step 2 in Section 3.2.2]{Ng-Complementary}, one can verify that $V_0 \in H^1(\Omega \setminus \partial B_{r_3})$ is a solution to 
\begin{equation*}
\left\{\begin{array}{c}\dive (s_0 A \nabla V_0) = f \mbox{ in } \Omega \setminus \partial B_{r_3},  \\[6pt]
[V_0] = g \quad \mbox{ and } \quad [A \nabla V_0 \cdot \eta] = h \mbox{ on } \partial B_{r_3}, \\[6pt]
V_0 = 0 \mbox{ on } \partial \Omega. 
\end{array}\right.
\end{equation*}
Set  
\begin{equation}\label{def-Wd-Thm2.3}
W_\delta = V_\delta - V_0 \mbox{ in } \Omega. 
\end{equation}
Then $W_\delta \in H^1_0(\Omega)$ is the unique solution to 
\begin{equation*}
\dive(s_\delta A \nabla W_\delta) =  - \dive\big( i \delta A \nabla V_0  1_{B_{r_2} \setminus B_{r_1} } \big)  \mbox{ in } \Omega.
\end{equation*}
Here and in what follows,  for a subset $D$ of $\mR^d$, $1_D$ denotes the characteristic function of $D$. 
Applying Lemma~\ref{lem-stability-Ng}, we have
\begin{equation}\label{haha3-Thm2.3}
\|W_\delta\|_{H^1(\Omega)} \le C \|V_0\|_{H^1(B_{r_2} \setminus B_{r_1})}. 
\end{equation}
The conclusion follows from \eqref{haha1-Thm2.3}, \eqref{haha2-Thm2.3}, \eqref{def-Wd-Thm2.3}, and \eqref{haha3-Thm2.3}. \proofend

\subsection{Proof of Theorem~\ref{thm3}}

\noindent \underline{Step 1:} Proof of \eqref{part1-thm3}. This is a consequence of Lemma~\ref{lem-compatible} with $g=h =0$. 

\medskip 
\noindent \underline{Step 2:} Proof of \eqref{part2-thm3}.  
 Without loss of generality, one might assume that $r_2 = 1$.  As in \cite{Ng-Complementary},  define 
\begin{equation*}
u_{1, \delta} = u_\delta \circ F^{-1} \mbox{ in } \mR^d \setminus B_{r_3}, 
\end{equation*}
and
\begin{equation*}
u_{2, \delta} = u_{1, \delta} \circ G^{-1} \mbox{ in } B_{r_3}. 
\end{equation*}
Let $\phi \in H^1_0(B_{r_3} \setminus B_{r_2})$ be the unique solution to 
\begin{equation}\label{def-phi-Thm2.3}
\Delta \phi = f \mbox{ in } B_{r_3} \setminus B_{r_2}, 
\end{equation}
and set 
\begin{equation*}
W = w - \phi \mbox{ in } B_{r_0} \setminus B_{r_2}. 
\end{equation*}
Then $W \in H^1(B_{r_0} \setminus B_{r_2})$ satisfies 
\begin{equation}\label{pro-V-S} 
\Delta W = 0 \mbox{ in } B_{r_0} \setminus B_{r_2}, \quad W = 0 \mbox{ on } \partial B_{r_2}, \quad \mbox{ and } \quad \partial_r W = - \partial_r \phi \mbox{ on } \partial B_{r_2}. 
\end{equation}
We now consider the case $d=2$ and $d=3$ separately. 

\medskip
\noindent \underline{Case 1:} $d=2$. 

\medskip 
Since $r_2 = 1$ and $W=0$ on $\partial B_{r_2}$, it follows that   
\begin{equation}\label{presentation-V-S}
W = g_{0} \ln r +  \sum_{\ell =1}^\infty \sum_{\pm} g_{\ell, \pm} (r^\ell - r^{-\ell}) e^{ \pm i \ell \theta} \mbox{ in } B_{r_0} \setminus B_{r_2}, 
\end{equation}
for some $g_0,  g_{\ell, \pm} \in  \mC$ ($\ell \ge 1$).  It is clear that, since $r_2 = 1 <  r_0$,  
\begin{equation}\label{finite-V}
\| W\|_{H^1(B_{r_0} \setminus B_{r_2})}^2 \sim |g_0|^2 + \sum_{\ell = 1}^\infty \sum_{\pm} \ell |g_{\ell, \pm} ^2| r_0^{2 \ell} < + \infty. 
\end{equation}
One of the keys  in the proof is the construction of  $W_\delta \in H^1(B_{r_3} \setminus B_{r_2})$ which is  defined as follows 
\begin{equation}\label{def-Vd-S}
W_\delta  = g_{0} \ln r + \sum_{\ell= 1}^{\infty} \sum_{\pm} \frac{g_{\ell, \pm}}{1 + \xi_\ell}  (r^\ell - r^{-\ell}) e^{\pm i \ell \theta}   \mbox{ in } B_{r_3} \setminus B_{r_2}, 
\end{equation}
where 
\begin{equation}\label{xil}
\xi_\ell  =  \delta^{1/2} (r_3 / r_0)^{\ell} \mbox{ for } \ell \ge 1. 
\end{equation}
Roughly speaking, $W_\delta$ is  the main part of the singularity of $u_\delta$.
From the definition of $W_\delta$, 
\begin{equation}\label{eq-Wdd}
\Delta W_\delta = 0 \mbox{ in } (B_{r_3} \setminus \bar B_{r_2}),  \quad W_\delta = 0 \mbox{ on } \partial B_{r_2}, 
\end{equation}
and 
\begin{equation}\label{vdelta-H1}
\|W_\delta \|_{H^1(B_{r_3} \setminus B_{r_2})}^2 \sim |g_0|^2 +  \sum_{\ell = 1}^\infty  \sum_{\pm} \frac{\ell|g_{\ell, \pm}|^2}{1 + \xi_\ell^2} r_3^{2 \ell}. 
\end{equation}
By \eqref{xil}, we have, if $\xi_\ell \le 1$, then
\begin{equation}\label{control1-v}
\frac{\ell |g_{\ell, \pm}|^2}{1 + \xi_\ell^2}  r_3^{2 \ell} \le 
 \ell   |g_{\ell, \pm}|^2r_3^{2\ell} \le \delta^{-1} \ell  |g_{\ell, \pm}|^2 r_0^{2\ell}, 
\end{equation}
and if $\xi_\ell \ge 1$, then
\begin{equation}\label{control2-v}
\frac{ \ell|g_{\ell, \pm}|^2  }{1 + \xi_\ell^2} r_3^{2 \ell} \le \ell |g_{\ell, \pm}|^2 r_3^{2 \ell} \xi_\ell^{-2} 
= \delta^{-1} \ell |g_{\ell, \pm}|^2 r_0^{2\ell}. 
\end{equation}
A combination of \eqref{finite-V},  \eqref{vdelta-H1}, \eqref{control1-v}, and \eqref{control2-v} yields  
\begin{equation}\label{vdelta-H1-1}
\|W_\delta \|_{H^1(B_{r_3} \setminus B_{r_2})} \le C \delta^{-1/2}. 
\end{equation}
Let $W_{1, \delta} \in H^1(\Omega)$ be the unique solution to 
\begin{equation*}\left\{
\begin{array}{c}\dive(s_\delta A \nabla W_{1, \delta}) = 0 \mbox{ in } \Omega \setminus \partial B_{r_2}, \\[6pt]
 [s_\delta A \nabla W_{1, \delta}\cdot \eta] = (-1 + i \delta) h_\delta  \mbox{ on } \partial B_{r_2}, \\[6pt]
 W_{1, \delta} = 0 \mbox{ on } \partial \Omega, 
\end{array}\right. 
\end{equation*}
where 
\begin{equation*}
h_\delta = - \partial_r (\phi + W_\delta)  \mbox{ on } \partial B_{r_2}, 
\end{equation*}
and let $W_{2, \delta} \in H^1(\Omega \setminus \partial B_{r_3})$ be the unique solution to 
\begin{equation*}\left\{
\begin{array}{c}
\dive(s_\delta A \nabla W_{2, \delta}) = f 1_{\Omega \setminus B_{r_3}} \mbox{ in } \Omega \setminus \partial B_{r_3}, \\[6pt]
[W_{2,\delta}] = \phi +  W_\delta \quad \mbox{ and } \quad  [A \nabla W_{2, \delta} \cdot \eta] = \partial_r \phi  +  \partial_r W_\delta   \mbox{ on } \partial B_{r_3},   \\[6pt]
W_{2, \delta} = 0 \mbox{ on } \partial \Omega.  
\end{array}\right.
\end{equation*}
Recall that,  for a subset $D$ of $\mR^d$, $1_D$ denotes the characteristic function of $D$. 
From \eqref{def-phi-Thm2.3}, \eqref{eq-Wdd}, and the fact $A=I$ in $B_{r_3} \setminus B_{r_2}$, we have 
\begin{equation}\label{u-V-delta}
u_\delta - ( \phi + W_\delta) {\bf 1}_{B_{r_3} \setminus B_{r_2}} = W_{1, \delta} + W_{2, \delta} \mbox{ in } \Omega. 
\end{equation}
Using \eqref{pro-V-S}, \eqref{presentation-V-S}, and \eqref{def-Vd-S}, we obtain 
\begin{equation*}
h_\delta = - \partial_r (\phi + W_\delta) = \partial_r (W- W_\delta)  = \partial_r \left( \sum_{\ell= 1}^{\infty} \sum_{\pm} \frac{\xi_\ell g_{\ell, \pm}}{1 + \xi_\ell}  (r^\ell - r^{-\ell}) e^{\pm i \ell \theta} \right) \mbox{ on } \partial B_{r_2}. 
\end{equation*}
Since $r_2 = 1$, it follows that 
\begin{equation}\label{hdelta-H12-S}
\| h_\delta \|_{H^{-1/2}(\partial B_{r_2})}^2 \sim  \sum_{\ell= 1}^{\infty} \sum_{\pm} \frac{\ell |\xi_\ell|^2 |g_{\ell, \pm}|^2}{1  + |\xi_\ell|^2}. 
\end{equation}
Using \eqref{xil}, we have, if $\xi_\ell \le 1$ then 
\begin{equation}\label{control1-h-S}
\frac{\ell |\xi_\ell|^2}{1 + |\xi_\ell|^2} |g_{\ell, \pm}|^2  \le 
\delta \ell  |g_{\ell, \pm}|^2 (r_3/ r_0)^{2 \ell} = \delta  \ell |g_{\ell, \pm}|^2 r_0^{2\ell} (r_3/ r_0^2)^{2 \ell} \le \delta \ell |g_{\ell, \pm}|^2 r_0^{2\ell},  
\end{equation}
since $r_0 >  \sqrt{r_2 r_3} = \sqrt{r_3}$,  and if $\xi_\ell \ge 1$  then 
\begin{equation}\label{control2-h-S}
\frac{ \ell |\xi_\ell|^2 }{1 + |\xi_\ell|^2} |g_{\ell, \pm}|^2 \le \ell |g_{\ell, \pm}|^2  = \ell  |g_{\ell, \pm}|^2 r_0^{2 \ell } r_0^{-2 \ell} \le  \delta \ell  |g_{\ell, \pm}|^2 r_0^{2 \ell }, 
\end{equation}
since $\delta^{1/2} r_0^{\ell} > \delta^{1/2} (r_3/ r_0)^{\ell} \ge 1$. 
A combination of  \eqref{hdelta-H12-S}, \eqref{control1-h-S}, and \eqref{control2-h-S} yields  
\begin{equation*}
\| h_\delta \|_{H^{-1/2}(\partial B_{r_2})} \le C \delta^{1/2} \| W\|_{H^{1/2}(\partial B_{r_0})} \le C \delta^{1/2}.
\end{equation*}
Applying Lemma~\ref{lem-stability-Ng}, we have
\begin{equation}\label{W1-1}
\|W_{1, \delta} \|_{H^1(\Omega)} \le (C/ \delta) \delta^{1/2} = C\delta^{-1/2}. 
\end{equation}
On the other hand,  from \eqref{vdelta-H1-1} and Lemma~\ref{lem-compatible}, we obtain
\begin{equation}\label{W2-1}
\|W_{2, \delta} \|_{H^1(B_{r_3} \setminus B_{r_3})} \le C \delta^{-1/2}. 
\end{equation}
The conclusion in the case $d=2$ now follows from \eqref{vdelta-H1-1}, \eqref{u-V-delta}, \eqref{W1-1}, and \eqref{W2-1}. 

\medskip
\noindent \underline{Case 2:} $d=3$.

\medskip 
Since $r_2 = 1$ and $W= 0$ on $\partial B_{r_2}$, it follows that   
\begin{equation*}
W = g_0 + \frac{\hat g_0}{r} +  \sum_{\ell = 1}^\infty \sum_{k = -\ell}^\ell g_{\ell, k} (r^\ell - r^{-\ell - 1}) Y_\ell^k(x/|x|) \mbox{ in } B_{r_0} \setminus B_{r_2}, 
\end{equation*}
for some $g_0, \hat g_0, g_{\ell, k} \in  \mC$. Here $Y_\ell^k$ is the spherical harmonic function of degree $\ell$ and of order $k$. Define $W_\delta \in H^1(B_{r_3} \setminus B_{r_2})$ as follows 
\begin{equation*}
W_\delta  =  g_0 + \frac{\hat g_0}{r} + \sum_{\ell= 1}^{\infty} \sum_{k = -\ell}^\ell \frac{g_{\ell, k}}{1 + \xi_\ell}  (r^\ell - r^{-\ell - 1}) Y_\ell^k(x/|x|)  \mbox{ in } B_{r_3} \setminus B_{r_2}, 
\end{equation*}
where 
\begin{equation*}
\xi_\ell  =  \delta^{1/2} (r_3 / r_0)^{\ell} \mbox{ for } \ell \ge 1. 
\end{equation*}
The proof now follows as in the two dimensional  case. The details are left to the reader. 
\proofend


\section{A connection between the blow up of the power and CALR. Proof of Theorem~\ref{thm1}} \label{sect-thm1}

We establish  a stronger result  than Theorem~\ref{thm1}: 
\begin{proposition}\label{pro1.1}
Let $d=2,  3$,  let $(\delta_n) \to 0$, $(g_n) \subset L^2(\Omega)$ 
with $\supp g_n \subset \Omega \setminus B_{r_2}$, and let $v_n \in H^1_0(\Omega)$ be the unique solution to 
\begin{equation*}
\dive (s_{\delta_n} A \nabla v_n) = g_n \mbox{ in } \Omega. 
\end{equation*}
Assume that  $s_0A$ is doubly complementary. Suppose that  $g_n \to g$ weakly in $L^2(\Omega)$, for some $g \in L^2(\Omega)$,  and
\begin{equation}\label{blowup-v-v}
\lim_{n \to \infty}\delta_n \| \nabla v_n\|_{L^2(B_{r_2} \setminus B_{r_1})} =  0.
\end{equation} 
Then $v_n \to v$ weakly in $H^1(\Omega \setminus B_{r_3})$ where $v \in H^1_0(\Omega)$ is the unique solution to 
\begin{equation*}
\dive (\hat A \nabla v) = g \mbox{ in } \Omega. 
\end{equation*}
\end{proposition}

\medskip
We first admit Proposition~\ref{pro1.1} and give 

\medskip
\noindent{\bf Proof of Theorem~\ref{thm1}.} Step 1: Proof of i). Since  $\delta_n \|\nabla v_{\delta_n} \|_{L^2(B_{r_2} \setminus B_{r_1})}^2 = 1$,  it follows from \eqref{PPP} that 
\begin{equation*}
\lim_{n \to + \infty} \delta_n \|\nabla v_{\delta_n} \|_{L^2(B_{r_2} \setminus B_{r_1})} = 0. 
\end{equation*}
On the other hand, since $\lim_{n \to + \infty} \delta_n \|\nabla u_{\delta_n} \|_{L^2(B_{r_2} \setminus B_{r_1})}^2 = + \infty$, we have
\begin{equation*}
\lim_{n \to + \infty} \| f_{\delta_n} \|_{L^2(\Omega)} = 0. 
\end{equation*}
The conclusion follows from Proposition~\ref{pro1.1}. 

\medskip
\noindent{Step 2:} Proof of ii). Since  $(\delta_n \|\nabla u_{\delta_n} \|_{L^2(B_{r_2} \setminus B_{r_1})}^2)$ is bounded,  it follows from \eqref{PPP} that 
\begin{equation*}
\lim_{n \to + \infty} \delta_n \|\nabla u_{\delta_n} \|_{L^2(B_{r_2} \setminus B_{r_1})} = 0. 
\end{equation*}
The conclusion now follows from Proposition~\ref{pro1.1}.  \proofend

\medskip
The rest of this section containing three subsections is devoted to the proof of Proposition~\ref{pro1.1}.  In the first subsection, we present the proof of Proposition~\ref{pro1.1} in the case $A = I$ in $B_{r_3} \setminus B_{r_2}$. This situation is already non-trivial since $A$ can be arbitrarily uniformly elliptic outside $B_{r_3}$; the standard separation of variables cannot be applied. Taking this simple but representative setting, we present the ideas of the proof of Proposition~\ref{pro1.1}. The proof essentially uses  the reflecting and removing localized singularity techniques introduced  in \cite{Ng-Complementary, Ng-Negative-Cloaking, Ng-superlensing}. The way to remove localized singularity in this context will lead us to develop the separation of variables  technique  for solving Cauchy problems in a general shell in Section~\ref{sect-separation}.  In Section~\ref{sect-pro1.1}, we give the proof of Proposition~\ref{pro1.1} in the form stated. To this end, we follow the strategy presented in Section~\ref{sect-motivation} and make use essentially the results in Section~\ref{sect-separation}. Due to the lack of the  orthogonality of plasmon modes, the analysis is more delicate.

\subsection{Proof of Proposition~\ref{pro1.1} in the case $A =I$ in $B_{r_3} \setminus B_{r_2}$} \label{sect-motivation}


Without loss of generality, one may assume that 
\begin{equation*}
r_3 = 1.
\end{equation*}
Using \eqref{PPP}, we derive from \eqref{blowup-v-v} that 
\begin{equation}\label{blowup-v}
\lim_{n \to \infty}\delta_n \|v_n\|_{H^1(\Omega)} =  0.
\end{equation} 
We now consider the case $d=2$ and $d=3$ separately.

\medskip
\noindent \underline{Case 1:} $d=2$. 

Define
\begin{equation*}
v_{1, n} = v_n \circ F^{-1} \mbox{ in } \mR^d \setminus B_{r_2}
\end{equation*}
and 
\begin{equation*}
v_{2, n} = v_{1, n} \circ G^{-1} \mbox{ in } B_{r_3}. 
\end{equation*}
It follows from  \eqref{def-DC}  and Lemma~\ref{lem-TO} that 
\begin{equation*}
\dive (A \nabla v_{1, n}) = \dive(A \nabla v_{2, n}) = 0 \mbox{ in } B_{r_3} \setminus B_{r_2}. 
\end{equation*}
Since
\begin{equation*}
A = I \mbox{ in } B_{r_3} \setminus B_{r_2},  
\end{equation*}
one can represent $v_{1, n}$ and $v_{2, n}$ in $B_{r_3} \setminus B_{r_2}$ as follows
\begin{equation}\label{Pmu1-cor}
v_{1, n} = c_0 + d_0 \ln r + \sum_{\ell = 1}^\infty \sum_{\pm} \big( c_{\ell, \pm} r^\ell + d_{\ell, \pm} r^{-\ell} \big) e^{\pm i \ell \theta}
\end{equation}
and 
\begin{equation}\label{Pmu2-cor}
v_{2, n} = e_0 + f_0 \ln r + \sum_{\ell = 1}^\infty \sum_{\pm}  \big(e_{\ell, \pm} r^\ell + f_{\ell, \pm} r^{-\ell} \big)e^{\pm i \ell \theta}, 
\end{equation}
for some $c_0, d_0, e_0, f_0, c_{\ell, \pm}, d_{\ell, \pm}, e_{\ell, \pm}, f_{\ell, \pm} \in \mC$ $(\ell \ge 1)$.
By Lemma~\ref{lem-TO}, we have
\begin{equation*}
v_{1, n} = v_{2, n} \quad \mbox{ and } \quad \partial_r  v_{1, n} = \frac{1}{1 - i \delta_n} \partial_r  v_{2, n}  \mbox{ on } \partial B_{r_3}. 
\end{equation*}
Since $r_3 = 1$, it follows that 
\begin{equation*}
c_{\ell, \pm} + d_{\ell, \pm} = e_{\ell, \pm} + f_{\ell, \pm} \quad \mbox{ and } c_{\ell, \pm} - d_{\ell, \pm} = \frac{1}{1 - i \delta_n} (e_{\ell, \pm} - f_{\ell, \pm} ) \quad \mbox{ for } \ell \ge 1, 
\end{equation*}
\begin{equation*}
c_0   = e_0  \quad \mbox{ and } \quad  d_0 = \frac{1}{1 - i \delta_n} f_0.  
\end{equation*}
This implies, for $\ell \ge 1$, 
\begin{equation*}
c_{\ell, \pm} = \frac{2 - i \delta_n}{2(1 - i \delta_n)} e_{\ell, \pm}  - \frac{i \delta_n}{2(1 - i \delta_n)} f_{\ell, \pm}  \quad  \mbox{ and } \quad d_{\ell, \pm} = \frac{2 - i \delta_n}{2(1 - i \delta_n)} f_{\ell, \pm}  - \frac{i \delta_n}{2(1 - i \delta_n)} e_{\ell, \pm}. 
\end{equation*}
We derive from \eqref{Pmu1-cor} and \eqref{Pmu2-cor} that 
\begin{equation}\label{v1-v2-Thm2.1-1}
v_{1, n} - v_{2, n} =   \frac{i \delta_n}{1 - i \delta_n} f_0 \ln r +  \frac{i \delta_n}{2(1 - i \delta_n)}  \sum_{\ell = 1}^\infty \sum_{\pm} (e_{\ell, \pm} - f_{\ell, \pm}) (r^{\ell} - r^{-\ell}) e^{\pm i \ell \theta}  \mbox{ in } B_{r_3} \setminus B_{r_2}.    
\end{equation}
It follows from \eqref{blowup-v} that 
\begin{equation*}
\lim_{n \to  \infty} \delta_n^2 \Big( \| v_{2, n}\|_{H^{1/2}(\partial B_{r_3})}^2 + \|\partial_r v_{2, n}\|_{H^{-1/2}(\partial B_{r_3})}^2 \Big) = 0
\end{equation*}
and
\begin{equation*}
\lim_{n \to  \infty} \delta_n^2 \Big( \| v_{2, n}\|_{H^{1/2}(\partial B_{r_2})}^2 + \|\partial_r v_{2, n}\|_{H^{-1/2}(\partial B_{r_2})}^2 \Big) = 0. 
\end{equation*}
Using \eqref{Pmu2-cor}, we obtain
\begin{equation}\label{bdr1-cor}
\lim_{n \to \infty } \delta_n^2 \Big( |e_0|^2  +\sum_{\ell = 1}^\infty \sum_{\pm} \ell |e_{\ell, \pm}|^2 r_3^{2 \ell} + |f_0|^2 + \sum_{\ell = 1}^\infty  \sum_{\pm} \ell|f_{\ell, \pm}|^2 r_3^{- 2 \ell} \Big)= 0 
\end{equation}
and 
\begin{equation}\label{bdr2-cor}
\lim_{n \to \infty } \delta_n^2 \Big( |e_0|^2   +\sum_{\ell = 1}^\infty \sum_{\pm} \ell|e_{\ell, \pm}|^2 r_2^{2 \ell} + |f_0|^2 + \sum_{\ell = 1}^\infty  \sum_{\pm} \ell|f_{\ell, \pm}|^2 r_2^{- 2 \ell} \Big)= 0. 
\end{equation}

We now use the removing localized singularity technique. Set 
\begin{equation}\label{vh-Thm2.1-1}
\hat v_n = - \frac{i \delta_n}{1 - i \delta_n}  f_0 \ln r -  \frac{i \delta_n}{2(1  - i \delta_n)}  \sum_{\ell =1}^\infty \sum_{\pm} (e_{\ell, \pm} - f_{\ell, \pm}) r^{-\ell} e^{ \pm i \ell \theta} \quad \mbox{ in } B_{r_3} \setminus B_{r_2}, 
\end{equation}
and  define $V_n$ in $\Omega$ as follows
\begin{equation}\label{def-Un-cor}
V_n = \left\{\begin{array}{cl} 
v_n & \mbox{ in } \Omega \setminus B_{r_3}, \\[6pt]
v_n - \hat v_n & \mbox{ in } B_{r_3}  \setminus B_{r_2}, \\[6pt]
v_{2, n} & \mbox{ in } B_{r_2}. 
\end{array}\right.
\end{equation}
Since $A = F_*A = G_*F_*A = I$ in $B_{r_3} \setminus B_{r_2}$,  we have, by Lemma~\ref{lem-TO}, 
\begin{equation}\label{eq-Un-cor}
\dive (\hat A \nabla V_n) = g_n \quad \mbox{ in } \Omega \setminus ( \partial B_{r_2} \cup \partial B_{r_3}), 
\end{equation}
where $\hat A$ is defined by \eqref{def-hatA}. 

We claim that 
\begin{equation}\label{cl1-cor}
\|[V_n] \|_{H^{1/2}(\partial B_{r_3})} + \|[\hat A \nabla V_n \cdot \eta] \|_{H^{-1/2}(\partial B_{r_3})} =  o(1)
\end{equation}
and 
\begin{equation}\label{cl2-cor}
\|[V_n] \|_{H^{1/2}(\partial B_{r_2})} + \|[\hat A \nabla V_n \cdot \eta] \|_{H^{-1/2}(\partial B_{r_2})} = o(1).
\end{equation}
Here and in what follows, $o(1)$ denotes a quantity converging to 0 as $n \to \infty$. 

\medskip 
We admit the claim and continue the proof. Combining  \eqref{eq-Un-cor}, \eqref{cl1-cor}, and \eqref{cl2-cor}  and using the fact that $V_n = 0$ on $\partial \Omega$ and $g_n \to g$ weakly in  $L^2(\Omega)$, we obtain 
\begin{equation*}
V_n \to v \mbox{ weakly in } H^1 \big(\Omega \setminus (\partial B_{r_3} \cup \partial B_{r_2}) \big),  
\end{equation*}
by the definition of $v$.  The conclusion follows since $v_n = V_n$ in $\Omega \setminus B_{r_3}$. 

\medskip
It remains to prove the claim. 

\medskip
\underline{Step 1.1}: Proof of \eqref{cl1-cor}.  Since $r_3  = 1$, we have, on $\partial B_{r_3}$,  
\begin{equation*}
[V_n] = \hat v_n =  - \sum_{\ell = 1}^\infty \sum_{\pm}\frac{i \delta_n}{2(1 - i \delta_n)} (e_{\ell, \pm} - f_{\ell, \pm}) r_3^{-\ell} e^{\pm i \ell \theta}. 
\end{equation*}
Since $r_3 = 1$, it follows from \eqref{bdr1-cor} and \eqref{bdr2-cor} that 
\begin{equation}\label{cl1-1-cor}
\| [V_n] \|_{H^{1/2}(\partial B_{r_3})} = o(1).  
\end{equation}
Similarly, 
\begin{equation}\label{cl1-2-cor}
\|[\hat A \nabla V_n \cdot \eta] \|_{H^{-1/2}(\partial B_{r_3})} =  o(1). 
\end{equation}
Claim \eqref{cl1-cor} is now a consequence of \eqref{cl1-1-cor} and \eqref{cl1-2-cor}. 

\medskip
\underline{Step 1.2}: Proof of \eqref{cl2-cor}. We have
\begin{equation*}
[V_n] = v_n - \hat v_n - v_{2, n} \mbox{ on } \partial B_{r_2}.  
\end{equation*}
This implies, since $v_n = v_{1, n}$ on $\partial B_{r_2}$, 
\begin{equation*}
[V_n]  =  v_{1, n}  - v_{2, n} - \hat v_n  \mbox{ on } \partial B_{r_2}. 
\end{equation*}
It follows from \eqref{v1-v2-Thm2.1-1} and \eqref{vh-Thm2.1-1} that 
\begin{equation*}
\|[V_n]\|_{H^{1/2}(\partial B_{r_2})} \le \Big\|  \frac{i \delta_n}{2(1 - i \delta_n)} \sum_{\ell = 1}^\infty \sum_{\pm} (e_{\ell, \pm} - f_{\ell, \pm})  r^\ell e^{ \pm i \ell \theta} \Big\|_{H^{1/2}(\partial B_{r_2})}. 
\end{equation*}
Since $r_3 =1$, we derive from \eqref{bdr1-cor} and \eqref{bdr2-cor}  that 
\begin{equation}\label{Un1-cor}
\|[V_n]\|_{H^{1/2}(\partial B_{r_2})} = o(1). 
\end{equation}
Similarly, using the fact that $\partial_r v_n = (1 - i \delta_n) \partial_r v_{1, n}$ and $\lim_{n \to \infty} \delta_n \| v_n\|_{H^1(\Omega)} = 0 $, we have
\begin{equation}\label{Un2-cor}
\|[\hat A \nabla V_n \cdot \eta ]\|_{H^{-1/2}(\partial B_{r_2})} = o(1). 
\end{equation}
A combination of \eqref{Un1-cor} and \eqref{Un2-cor} yields \eqref{cl2-cor}.

\medskip
\noindent \underline{Case 2:} $d=3$. 

The proof in the three dimensional case follows similarly as the one in the two dimensional case. We just note here that, in three dimensions, $v_{1,n}$ and $v_{2, n}$ can be represented as follows
\begin{equation*}
v_{1, n} = c_{0, 0} + \frac{d_{0, 0}}{r} +  \sum_{\ell=1}^\infty \sum_{k=-\ell}^\ell (c_{\ell, k} r^\ell + d_{\ell, k} r^{-\ell-1}) Y^k_\ell(x/|x|) \quad \mbox{ in } B_{r_3} \setminus B_{r_2}
\end{equation*}
and 
\begin{equation*}
v_{2, n} =  e_{0, 0} + \frac{f_{0, 0}}{r} +  \sum_{\ell=1}^\infty \sum_{k=-\ell}^\ell(e_{\ell, k} r^\ell + f_{\ell, k} r^{-\ell-1}) Y^k_\ell(x/|x|) \quad \mbox{ in } B_{r_3} \setminus B_{r_2}, 
\end{equation*}
for some $c_{\ell, k}, \; d_{\ell, k}, \; e_{\ell, k}, \; f_{\ell, k} \in \mC$. 
\proofend

\subsection{Separation of variables approach for Cauchy problems in a general shell} \label{sect-separation}

In this section, we state variants of \eqref{Pmu1-cor} and \eqref{Pmu2-cor} for a general core-shell structure, i.e., $A$ is not required to be $I$ in $B_{r_3} \setminus B_{r_2}$. Using these variants, we can extend the method used in Section~\ref{sect-motivation} for a general core-shell structure in Section~\ref{sect-pro1.1}. We have

\begin{proposition}\label{pro1}
Let $d = 2, \, 3$,  $0 < R_1 < R_2$, and let $a \in [C^3(\overline{B_{R_2} \setminus B_{R_1}})]^{d \times d}$ be symmetric and uniformly elliptic. Set $R_3 = R_2^2/ R_1$ and let $K: B_{R_2} \setminus B_{R_1} \to B_{R_3} \setminus B_{R_2}$ be the Kelvin transform with respect to  $\partial B_{R_2}$, i.e., $K(x) = x R_2^2/ |x|^2$. Define
\begin{equation}\label{def-a1}
a_1 = \left\{\begin{array}{cl} a &  \mbox{ in } B_{R_2} \setminus B_{R_1}, \\[6pt]
K_* a & \mbox{ in } B_{R_3} \setminus B_{R_2}, \\[6pt]
I & \mbox{ in } B_{R_1}. 
\end{array}\right.
\end{equation}
Let $v_\ell \in H^1(B_{R_3})$ ($\ell \ge 1$) be a solution to  
\begin{equation*}
\dive (a_1 \nabla v_\ell ) = 0 \mbox{ in } B_{R_3}, 
\end{equation*}
and set $v_0 = 1$ in  $B_{R_3}$.  Let $w_\ell \in H^1(B_{R_2} \setminus B_{R_1})$ ($\ell \ge 1$)   be  the reflection of $v_\ell$ through $\partial B_{R_2}$ by $K^{-1}$, i.e., 
\begin{equation*}
w_\ell = v_\ell \circ K \mbox{ in } B_{R_2} \setminus B_{R_1}, 
\end{equation*}
and denote $w_0 \in H^1(B_{R_3} \setminus B_{R_2})$ the unique solution to  
\begin{equation*}
\dive (a \nabla w_0) = 0 \mbox{ in } B_{R_2} \setminus B_{R_1}, \quad w_{0} = 1 \mbox{ on } \partial B_{R_2}, \quad \mbox{ and } \quad w_0 = 0 \mbox{ on } \partial B_{R_1}.   
\end{equation*}
Then, for $\ell \ge 1$, 
\begin{equation}\label{pro-wl-1}
\dive (a \nabla w_\ell) = \dive (a \nabla v_\ell) = 0 \mbox{ in } B_{R_2} \setminus B_{R_1},
\end{equation}
\begin{equation}\label{pro-wl-2}
w_\ell =  v_\ell \mbox{ on } \partial B_{R_2}, \quad  \mbox{ and } \quad  a \nabla w_\ell \cdot \frac{ x}{ |x|} = - a \nabla v_\ell \cdot \frac{x}{|x|} \mbox{ on } \partial B_{R_2}. 
\end{equation}
Assume that $\{v_\ell \}_0^\infty$ is  dense in $H^{1/2}(\partial B_{R_3})$. We have, with respect to  $H^1(B_{R_2} \setminus B_{R_1})$-norm,  

\begin{enumerate}

\item[1)]
\begin{equation*}
\Big\{v_\ell - w_\ell ; \; \ell \ge 0  \Big\} \mbox{ is dense in }  \Big\{v \in H^1(B_{R_2} \setminus B_{R_1}); \; \dive(a \nabla v) = 0 \mbox{ and }  v = 0 \mbox{ on } \partial B_{R_2} \Big\}.
\end{equation*}

\item[2)]  
\begin{multline*}
\Big\{1 \Big\} \cup \Big\{v_\ell + w_\ell ; \; \ell \ge 1  \Big\} \mbox{ is dense in } \\[6pt] \Big\{v \in H^1(B_{R_2} \setminus B_{R_1});  \; \dive(a \nabla v) = 0 \mbox{ and }  a \nabla v \cdot \eta = 0 \mbox{ on } \partial B_{R_2} \Big\}.
\end{multline*}

\item[3)]
\begin{equation*}
\Big\{v_\ell, w_\ell ; \; \ell \ge 0  \Big\} \mbox{ is dense in } \Big\{v \in H^1(B_{R_2} \setminus B_{R_1}) \; \dive(a \nabla v) = 0 \Big\}.  
\end{equation*}
\end{enumerate}

\end{proposition}

The proof of Proposition~\ref{pro1} is given in the appendix. 

\medskip

The existence of $v_{\ell}$ and $w_\ell$, their density properties, and  \eqref{pro-wl-1} and \eqref{pro-wl-2} can be considered as the existence of surface plasmons for complementary media, a fact which  can be used elsewhere; see e.g.,  \cite{Grieser0, Grieser, Maier} for discussions on surface plasmons and their applications.  The choice of $a_1$ is to ensure such properties. 


\subsection{Proof of Proposition~\ref{pro1.1}}\label{sect-pro1.1}

Using \eqref{PPP}, we derive from \eqref{blowup-v-v} that 
\begin{equation}\label{blowup-v-1}
\lim_{n \to \infty}\delta_n \|v_n\|_{H^1(\Omega)} =  0.
\end{equation} 
Define 
\begin{equation*}
v_{1, n} = v_n \circ F^{-1} \mbox{ in } B_{r_4} \setminus B_{r_3}
\end{equation*}
and 
\begin{equation*}
v_{2, n} = v_{1, n} \circ G^{-1} \mbox{ in } B_{r_3}. 
\end{equation*}
Using \eqref{def-DC} and applying Lemma~\ref{lem-TO}, we obtain
\begin{equation*}
\dive (A \nabla v_{1, n}) = \dive(A \nabla v_{2, n}) = 0 \mbox{ in } B_{r_3} \setminus B_{r_2}
\end{equation*}
and
\begin{equation*}
v_{1, n} = v_{2, n} \quad \mbox{ and } \quad A \nabla v_{1, n} \cdot \eta = \frac{1}{1 - i \delta_n} A \nabla v_{2, n} \cdot \eta \mbox{ on } \partial B_{r_3}. 
\end{equation*}
Set 
\begin{equation*}
\hat r = r_3^2/ r_2
\end{equation*}
and let $K: B_{r_3} \setminus B_{r_2} \to B_{\hat r} \setminus B_{r_3}$ be the Kelvin transform with respect to $\partial B_{r_3}$. Define
\begin{equation}\label{def-a1-1-1}
A_1 = \left\{\begin{array}{cl} A &  \mbox{ in } B_{r_3} \setminus B_{r_2}, \\[6pt]
{K}_* A & \mbox{ in } B_{\hat r} \setminus B_{r_3}, \\[6pt]
I & \mbox{ in } B_{r_2}. 
\end{array}\right.
\end{equation}
Let $v_\ell \in H^1(B_{\hat r})$ ($\ell \ge 1$) be a solution to  
\begin{equation*}
\dive (A_1 \nabla v_\ell ) = 0 \mbox{ in } B_{\hat r}, 
\end{equation*}
and set $v_0 = 1$ in  $B_{\hat r}$.  Define $w_\ell \in H^1(B_{r_3} \setminus B_{r_2})$ ($\ell \ge 1$)  the reflection of $v_\ell$ through $\partial B_{r_3}$ by $K^{-1}$, i.e., 
\begin{equation}\label{def-www}
w_\ell = v_\ell \circ K \mbox{ in } B_{r_3} \setminus B_{r_2}, 
\end{equation}
and denote $w_0 \in H^1(B_{r_3} \setminus B_{r_2})$ the unique solution to  
\begin{equation*}
\dive (A \nabla w_0) = 0 \mbox{ in } B_{r_3} \setminus B_{r_2}, \quad w_{0} = 1 \mbox{ on } \partial B_{r_3}, \quad \mbox{ and } \quad w_0 = 0 \mbox{ on } \partial B_{r_2}.   
\end{equation*}
We assume in additional that $\{ v_\ell \}_0^\infty$ is  an orthogonal basis  of $H^{1/2}(\partial B_{\hat r})$.  In particular, we have
\begin{equation}\label{ortho}
\int_{\partial B_{\hat r}} v_{\ell} = 0 \mbox{ for } \ell \ge 1 \;  \footnote{In the case $d=2$ and $r_3 =1$, $v_\ell$ and $w_\ell$ can be seen as a replacement of  $r^\ell e^{\pm i \ell \theta}$ and $r^{-\ell} e^{\pm i \ell \theta}$ respectively.}. 
\end{equation}
For $m \ge 0$, let $P_m$ be the projection from $H^1(B_{r_3} \setminus B_{r_2})$ into  the $\mbox{span} \big\{v_\ell, \, w_\ell; \; 0 \le  \ell \le m \big\}$ with respect to $H^1(B_{r_3} \setminus B_{r_2})$-norm. 
By  Proposition~\ref{pro1},  there exists $m$  such that 
\begin{equation}\label{projection}
\|v_{1, n} - P_m v_{1, n} \|_{H^1(B_{r_3} \setminus B_{r_2})} + \|v_{2, n} - P_m v_{2, n} \|_{H^1(B_{r_3} \setminus B_{r_2})} \le \delta_n^2. 
\end{equation}
We have, in $B_{r_3} \setminus B_{r_2}$, 
\begin{equation}\label{Pmu1}
P_m v_{1, n} = \sum_{\ell = 0}^{m} \big( c_{\ell} v_\ell + d_\ell w_\ell \big)
\end{equation}
and 
\begin{equation}\label{Pmu2}
P_m v_{2, n} = \sum_{\ell = 0}^{ m} \big(e_{\ell} v_\ell + f_\ell w_\ell \big), 
\end{equation}
for some $c_\ell, d_\ell, e_\ell, f_\ell \in \mC$ $(0 \le \ell \le m)$.
Define $(D_\ell)_0^m, (N_\ell)_0^m \subset \mC$ as follows 
\begin{equation}\label{dd1}
c_\ell + d_\ell = e_\ell + f_\ell + D_\ell \quad \mbox{ and } c_\ell - d_\ell = \frac{1}{1 - i \delta_n} (e_\ell - f_\ell ) + N_\ell \quad \mbox{ for } 1\le  \ell \le m
\end{equation}
and 
\begin{equation}\label{dd2}
c_0  + d_0= e_0 + f_0  + D_0\quad \mbox{ and } d_0 = \frac{1}{1 - i \delta_n} f_0 + N_0.  
\end{equation}
It follows from \eqref{pro-wl-2} that 
\begin{equation}\label{def-Dl}
P_m v_{1, n} - P_m v_{2, n} = \sum_{\ell =0}^{ m} D_\ell v_\ell \mbox{ on } \partial B_{r_3}
\end{equation}
and 
\begin{equation}\label{defNl}
a \nabla P_m v_{1, n} \cdot \eta -  \frac{1}{1 - i \delta_n} a \nabla  P_m v_{2, n}  \cdot \eta = N_0 \, a \nabla w_0   \cdot \eta + \sum_{\ell = 1}^{m} N_\ell \, a \nabla v_\ell \cdot \eta \mbox{ on } \partial B_{r_3}. 
\end{equation}
From \eqref{dd1} and \eqref{dd2}, we have, for $1 \le \ell \le m$, 
\begin{equation*}
c_\ell = \frac{2 - i \delta_n}{2(1 - i \delta_n)} e_\ell  - \frac{i \delta_n}{2(1 - i \delta_n)} f_\ell + \frac{D_\ell + N_\ell}{2}, \quad d_\ell = \frac{2 - i \delta_n}{2(1 - i \delta_n)} f_\ell  - \frac{i \delta_n}{2(1 - i \delta_n)} e_\ell + \frac{D_\ell - N_\ell}{2},  
\end{equation*}
and
\begin{equation*}
c_0 =  e_0  - \frac{ i \delta_n}{1 - i \delta_n} f_0 + D_0 - N_0 \quad  \mbox{ and } \quad d_0 = \frac{1}{1 - i \delta_n} f_0 + N_0.
\end{equation*}
We derive from \eqref{Pmu1} and \eqref{Pmu2} that  
\begin{align}\label{P1-P2}
P_m v_{1, n} - P_m v_{2, n} = & \frac{i \delta_n}{2(1 - i \delta_n)}  \sum_{\ell =1}^m (e_\ell - f_\ell) (v_\ell - w_\ell) + \sum_{\ell=1}^m \left( \frac{D_\ell + N_\ell}{2} v_\ell + \frac{D_\ell - N_\ell}{2} w_\ell \right) \nonumber \\[6pt]
&  + \Big( - \frac{i \delta_n}{1 - i \delta_n} f_0 + D_0 - N_0\Big) + \Big( \frac{i \delta_n}{1 - i \delta_n} f_0 + N_0 \Big) w_0.   
\end{align}
From \eqref{blowup-v-1} and \eqref{projection}, we have
\begin{equation}\label{norm-projection}
\| P_m v_{2, n} \|_{H^{1/2}(\partial B_{r_3})} =  \delta_n^{-1} o(1) \quad \mbox{ and } \quad \| P_m v_{2, n} \|_{H^{1/2}(\partial B_{r_2})} =  \delta_n^{-1} o(1).   
\end{equation}
Here and in what follows in this proof, $o(1)$ denotes a quantity converging to 0 as $n \to \infty$. 
Since $v_\ell = w_\ell$  on $\partial B_{r_3}$ for $\ell \ge 1$,  it follows from \eqref{Pmu2} and \eqref{norm-projection} that 
\begin{equation}\label{est1-1}
\Big\| \sum_{\ell =0}^m (e_\ell + f_\ell)v_\ell \Big\|_{H^{1/2}(\partial B_{r_3})} = \delta_n^{-1} o(1)
\end{equation}
and
\begin{equation}\label{est1-2}
\Big\| \sum_{\ell=0}^m (e_\ell v_\ell + f_\ell w_\ell) \Big\|_{H^{1/2}(\partial B_{r_2})}  = \Big\| \sum_{\ell=0}^m (e_\ell + f_\ell)v_\ell + \sum_{\ell=0}^m f_\ell(w_\ell - v_\ell) \Big\|_{ H^{1/2}(\partial B_{r_2})} = \delta_n^{-1} o(1). 
\end{equation}
Since, for $\ell \ge 0$,  
\begin{equation*}
\dive (A_1 \nabla v_\ell )  = 0 \mbox{ in } B_{r_3}, 
\end{equation*}
we have
\begin{equation}\label{haha1}
\Big\| \sum_{\ell =0}^m (e_\ell + f_\ell)v_\ell \Big\|_{H^{1/2}(\partial B_{r_2})} \le C \Big\| \sum_{\ell =0}^m (e_\ell + f_\ell)v_\ell \Big\|_{H^{1/2}(\partial B_{r_3})}. 
\end{equation}
Here and in what follows in this proof, $C$ denotes a positive constant independent of $\delta_n$, $u_n$, $g_n$, and $\ell$.  A combination of \eqref{est1-1}, \eqref{est1-2}, and \eqref{haha1} yields
\begin{equation}\label{est1-3}
\Big\| \sum_{\ell =0}^{ m} f_\ell (w_\ell - v_\ell) \Big\|_{ H^{1/2}(\partial B_{r_2})} =  \delta_n^{-1}o(1). 
\end{equation}
Using \eqref{ortho} and applying  Lemma~\ref{lem-stability} below with  $v = -\sum_{\ell \ge 1}^m f_\ell v_\ell$, $c = f_0$, $R_1 = r_2$, and $R_2 = r_3$, we deduce from \eqref{est1-3} that 
\begin{equation}\label{est1-4}
|f_0| +  \Big\| \sum_{\ell =1}^m f_\ell v_\ell \Big\|_{H^{1/2}(\partial B_{r_3})} + \Big\| \sum_{\ell =0}^{ m} f_\ell w_\ell  \Big\|_{ H^{1/2}(\partial B_{r_2})} =  \delta_n^{-1}o(1). 
\end{equation}
We also use here the fact that $w_0 = 0$ on $\partial B_{r_2}$. 
This implies, by  \eqref{est1-1},  
\begin{equation}\label{est1-5}
 \Big\| \sum_{\ell=0}^m e_\ell  v_\ell \Big\|_{H^{1/2}(\partial B_{r_3})}  = \delta_n^{-1} o(1).
\end{equation}
From \eqref{est1-4} and \eqref{est1-5}, we obtain
\begin{equation}\label{est1-6}
\Big\| \sum_{\ell=0}^m e_\ell  v_\ell \Big\|_{H^{1/2}(\partial B_{r_3})}  + |f_0| + \Big\| \sum_{\ell =0}^{ m} f_\ell w_\ell  \Big\|_{ H^{1/2}(\partial B_{r_2})}  = \delta_n^{-1} o(1).
\end{equation}
Since $\dive (A_1 \nabla v_\ell) = 0$ in $B_{\hat r}$ for $\ell \ge 1$, $v_0 = 1$,  and $A_1 = A$ in $B_{r_3} \setminus B_{r_2}$, 
\begin{equation}\label{est1-7}
\Big\| \sum_{\ell =1}^m e_\ell  A \nabla v_\ell \cdot \eta \Big\|_{H^{-1/2}(\partial B_{r_3})}  \le C \Big\| \sum_{\ell=0}^m e_\ell  v_\ell \Big\|_{H^{1/2}(\partial B_{r_3})}.
\end{equation}
From \eqref{def-a1-1-1} and  \eqref{def-www}, we have
\begin{multline}\label{est1-8}
\Big\| \sum_{\ell=1}^m f_\ell A \nabla w_\ell  \cdot \eta \Big\|_{H^{-1/2}(\partial B_{r_2})} \le C 
\Big\| \sum_{\ell=1}^m f_\ell A_1 \nabla v_\ell  \cdot \eta \Big\|_{H^{-1/2}(\partial B_{\hat r})}  \\[6pt]
 \le C 
\Big\| \sum_{\ell=1}^m f_\ell  v_\ell   \Big\|_{H^{1/2}(\partial B_{\hat r})} 
  \le C  \Big\| \sum_{\ell =1}^{ m} f_\ell w_\ell  \Big\|_{ H^{1/2}(\partial B_{r_2})}. 
\end{multline}
Recall that $w_0 = 0$ on $\partial B_{r_2}$. A combination of \eqref{est1-6}, \eqref{est1-7}, and \eqref{est1-8} yields 
\begin{equation}\label{est1-4-1}
\Big\| \sum_{\ell =1}^m e_\ell  A \nabla v_\ell \cdot \eta \Big\|_{H^{-1/2}(\partial B_{r_3})}  +  \Big\| \sum_{\ell=1}^m f_\ell A \nabla w_\ell  \cdot \eta \Big\|_{H^{-1/2}(\partial B_{r_2})} = \delta_n^{-1} o(1). 
\end{equation}

We are ready to use the removing localized singularity. Set, in $B_{r_3} \setminus B_{r_2}$,  
\begin{align*}
\hat v_n = & - \sum_{\ell = 1}^m \frac{i \delta_n}{2(1 - i \delta_n)} (e_\ell - f_\ell) w_\ell  + \sum_{\ell=1}^m \left( \frac{D_\ell + N_\ell}{2} v_\ell + \frac{D_\ell - N_\ell}{2} w_\ell \right) \\[6pt]
& +  \Big(- \frac{i \delta_n}{1 - i \delta_n} f_0 + D_0 - N_0 \Big) + \Big(\frac{i \delta_n}{1 - i \delta_n} f_0 + N_0 \Big) w_0 - \frac{i \delta_n}{2(1 - i \delta_n)} (e_0 - f_0) v_0.  
\end{align*}
It follows from \eqref{P1-P2} that 
\begin{equation}\label{pro-diff}
P_m v_{1, n} - P_m v_{2, n}  =   \frac{i \delta_n}{2(1 - i \delta_n)} \sum_{\ell = 0}^m (e_\ell - f_\ell)  v_\ell   + \hat v_n \quad \mbox{ in } B_{r_3} \setminus B_{r_2}.  
\end{equation}
Define 
\begin{equation}\label{def-Un}
V_n = \left\{\begin{array}{cl} 
v_n & \mbox{ in } \Omega \setminus B_{r_3}, \\[6pt]
v_n - \hat v_n & \mbox{ in } B_{r_3}  \setminus B_{r_2}, \\[6pt]
v_{2, n} & \mbox{ in } B_{r_2}. 
\end{array}\right.
\end{equation}
We have 
\begin{equation}\label{eq-Un}
\dive (\hat A \nabla V_n) = g_n \mbox{ in } \Omega \setminus ( \partial B_{r_2} \cup \partial B_{r_3}). 
\end{equation}

We claim that 
\begin{equation}\label{cl1}
\|[V_n] \|_{H^{1/2}(\partial B_{r_3})} + \|[\hat A \nabla V_n \cdot \eta] \|_{H^{1/2}(\partial B_{r_3})} =  o(1)
\end{equation}
and 
\begin{equation}\label{cl2}
\|[V_n] \|_{H^{1/2}(\partial B_{r_2})} + \|[\hat A \nabla V_n \cdot \eta] \|_{H^{1/2}(\partial B_{r_2})} = o(1).
\end{equation}
Admitting \eqref{cl1} and \eqref{cl2}, we derive  that $V_n \to v$ weakly in $H^1(\Omega \setminus (\partial B_{r_2} \cup \partial B_{r_3}))$ as in Section~\ref{sect-motivation}. The conclusion now follows from \eqref{def-Un}. 

\medskip It remains to prove \eqref{cl1} and \eqref{cl2}. 

\medskip 

\noindent \underline{Step 1}: Proof of \eqref{cl1}.  We have, on $\partial B_{r_3}$,  
\begin{equation*}
[V_n] = \hat v_n = - \sum_{\ell = 0}^m \frac{i \delta_n}{2(1 - \delta_n)} (e_\ell - f_\ell) v_\ell  + \sum_{\ell =0}^m D_\ell v_\ell.  
\end{equation*}
Here we use the fact that $w_\ell = v_\ell$ ($\ell \ge 0$) on $\partial B_{r_3}$.  We derive from \eqref{projection},  \eqref{def-Dl},  \eqref{est1-4},  and  \eqref{est1-5} that
\begin{equation}\label{cl1-1}
\|[V_n] \|_{H^{1/2}(\partial B_{r_3})} = o(1). 
\end{equation}
Similarly, using the fact that $A \nabla v_\ell \cdot \eta = - A \nabla w_{\ell} \cdot \eta$ on $\partial B_{r_3}$ for $\ell \ge 1$, and $f_0 =\delta_n^{-1} o(1)$,  
we derive from \eqref{projection}, \eqref{defNl}, and \eqref{est1-4-1} that 
\begin{equation}\label{cl1-2}
\|[A \nabla V_n \cdot \eta] \|_{H^{1/2}(\partial B_{r_3})} =  o(1). 
\end{equation}
A combination of \eqref{cl1-1} and \eqref{cl1-2} yields \eqref{cl1}. 

\medskip
\noindent \underline{Step 2}: Proof of \eqref{cl2}. We have, on $\partial B_{r_2}$, 
\begin{equation*}
[V_n] = v_n - \hat v_n - v_{2, n}. 
\end{equation*}
It follows that, on $\partial B_{r_2}$
\begin{equation*}
[V_n]  = v_n - v_{1, n} + v_{1, n} - P_m v_{1, n} + P_m v_{1, n} - P_m v_{2, n} + P_m v_{2, n} - v_{2, n} - \hat v_n. 
\end{equation*}
Since $v_n = v_{1, n}$ on $\partial B_{r_2}$, we derive from \eqref{projection} and \eqref{pro-diff} that 
\begin{equation*}
\|[V_n]\|_{H^{1/2}(\partial B_{r_2})} \le \delta_n^2 + \Big\|  \frac{i \delta_n}{2(1 - i \delta_n)} \sum_{\ell = 0}^m (e_\ell - f_\ell)  v_\ell   \Big\|_{H^{1/2}(\partial B_{r_2})}
\end{equation*}
From  \eqref{est1-4} and \eqref{est1-6}, we obtain 
\begin{equation}\label{Un1}
\|[V_n]\|_{H^{1/2}(\partial B_{r_2})} = o(1). 
\end{equation}
Similarly,
\begin{equation}\label{Un2}
\|[\hat A \nabla V_n \cdot \eta ]\|_{H^{-1/2}(\partial B_{r_2})} = o(1). 
\end{equation}
A combination of \eqref{Un1} and \eqref{Un2} yields \eqref{cl2}.

\medskip 
The proof is complete. 
\proofend

\medskip
In the proof of Proposition~\ref{pro1.1}, we used the following lemma. 

\begin{lemma}\label{lem-stability} Let $d=2, 3$, $0 < R_1 < R_2$, and let $a$ be a uniformly elliptic matrix-valued function defined in $B_{R_2} \setminus B_{R_1}$. Set $R_3 = R_2^2/ R_1$ and let $K: B_{R_2} \setminus B_{R_1} \to B_{R_3} \setminus B_{R_2}$ be the Kelvin transform with respect to  $\partial B_{R_2}$. Define
\begin{equation*}
a_1 = \left\{\begin{array}{cl} a &  \mbox{ in } B_{R_2} \setminus B_{R_1}, \\[6pt]
K_* a & \mbox{ in } B_{R_3} \setminus B_{R_2}, \\[6pt]
I & \mbox{ in } B_{R_1}. 
\end{array}\right.
\end{equation*}
Let 
$v \in H^1(B_{R_3})$ be such that 
\begin{equation*}
 \int_{\partial B_{R_3}} v = 0 \quad \mbox{ and } \quad \dive(a_1 \nabla v) = 0 \mbox{ in } B_{R_3}, 
\end{equation*}
and let $w \in H^1(B_{R_2} \setminus B_{R_1})$ be the reflection of $w$ by $K^{-1}$ through $\partial B_{R_2}$, i.e., 
\begin{equation*}
w = v \circ K \mbox{ in } B_{R_2} \setminus B_{R_1}. 
\end{equation*}
We have, for all $c \in \mC$, 
\begin{equation*}
\|v \|_{ H^{1/2}(\partial B_{R_2})} +  |c| \le  C \| v - w  + c\|_{H^{1/2}(\partial B_{R_1})}, 
\end{equation*}
where $C$ is a  positive constant independent of $v$ and $c$. 
\end{lemma}

\noindent{\bf Proof.} We prove Lemma~\ref{lem-stability} by contradiction. Assume that the conclusion is not true. Then there exists a sequence $(v_n) \subset H^1(B_{R_3})$ and $(c_n) \subset \mC$ such that 
\begin{equation}\label{vn-11}
\dive(a_1 \nabla v_n) = 0 \mbox{ in } B_{R_3}, 
\end{equation}
\begin{equation}\label{vn-22}
 \int_{\partial B_{R_3}} v_n = 0, \quad \|v_n \|_{H^{1/2}(\partial B_{R_2})} + |c_n| = 1, \quad \mbox{ and } \quad \lim_{n \to \infty}\| v_n - w_n + c_n \|_{H^{1/2}(\partial B_{R_1})} =0. 
\end{equation}
Here $w_n$ is the reflection of $v_n$ with respect to  $\partial B_{R_2}$ by $K^{-1}$.  From \eqref{vn-22}, we have
\begin{equation*}
\| v_n + c_n \|_{H^{1/2}(\partial B_{R_1})} \le C. 
\end{equation*}
In this proof, $C$ denotes a positive constant independent of $n$. It follows from \eqref{vn-22} that 
\begin{equation*}
\|w_n \|_{H^{1/2}(\partial B_{R_1})} \le C; 
\end{equation*}
which implies, by the definition of $w_n$, 
\begin{equation*}
\| v_n\|_{H^{1/2}(\partial B_{R_3})} \le C. 
\end{equation*}
Without loss of generality, one might assume that $v_n \to v$ weakly in $H^1(B_{R_3})$, $v_n \to v$ in $H^1_{\loc} (B_{R_3})$, and $c_n \to c \in \mC$. Moreover, from \eqref{vn-11} and \eqref{vn-22}, we have
\begin{equation}\label{v-11}
\dive(a_1 \nabla v) = 0 \mbox{ in } B_{R_3}, 
\end{equation}
\begin{equation}\label{v-22}
\int_{\partial B_{R_3}} v = 0,  \quad \mbox{ and } \quad \|v \|_{H^{1/2}(\partial B_{R_2})} = 1. 
\end{equation}
Let $w$ be the reflection of $v$ with respect to  $\partial B_{R_2}$ by $K^{-1}$. Since $v_n \to v$ in $H^1(B_{R_2})$, it follows from \eqref{vn-22} that  
\begin{equation*}
\lim_{n \to \infty} \|w_n - w\|_{H^{1/2}(\partial B_{R_1})} = 0; 
\end{equation*}
which  implies 
\begin{equation*}
\lim_{n \to \infty} \|v_n - v\|_{H^{1/2}(\partial B_{R_3})} = 0.  
\end{equation*}
From \eqref{vn-22}, we have 
\begin{equation*}
v - w + c = 0 \mbox{ on } \partial B_{R_1}. 
\end{equation*}
It follows from Lemma~\ref{lem1} in the appendix that  $v =0$ and $c =0$.  Here we use the fact that $\dsp \int_{\partial B_{R_3}}  v=0$. This contradicts \eqref{v-22}.  \proofend

\section{Cloaking a source via anomalous localized resonance}\label{sect-cloaking}

In this section, we describe how to use the theory CALR discussed previously to cloak a source $f$ concentrating on an arbitrary bounded smooth manifold of codimension 1 in an arbitrary medium.   
Without loss of generality, one may assume that the medium is contained in $B_{r_3} \setminus B_{r_2}$ and characterized by a matrix $a$ which is assumed smooth and uniformly elliptic in $\overline{B_{r_3} \setminus B_{r_2}}$ for some $0< r_2 <  r_3$.  Assume that $f$ concentrates on $\partial D$ for some bounded smooth open subset $D \subset \subset B_{r_3} \setminus B_{r_2}$. One might assume as well that $D \subset \subset B_{r_*}$ where $r_*$ is the constant coming from Theorem~\ref{thm2} since one can choose $r_3$ large enough (see \cite[Lemma 1]{Ng-Negative-Cloaking}). 
Define $r_1 = r_2^2/ r_3$.  Let $F: B_{r_2} \setminus \{0 \} \to \mR^d \setminus B_{r_2}$ and $G: \mR^d \setminus B_{r_3} \to B_{r_3} \setminus \{0 \}$ be the Kelvin transform with respect to $\partial B_{r_2}$ and $\partial B_{r_3}$ respectively. Note that $G \circ F (x) = (r_2^2 / r_1^2) x $. Define 
\begin{equation}\label{A-cloak}
A = \left\{\begin{array}{cl}  a & \mbox{ in } B_{r_3} \setminus B_{r_2}, \\[6pt]
F^{-1}_*a  & \mbox{ in } B_{r_2} \setminus B_{r_1}, \\[6pt]
F^{-1}_* G^{-1}_*a  & \mbox{ in } B_{r_1} \setminus B_{r_1^2/ r_2}, \\[6pt]
I & \mbox{ otherwise}.  
\end{array}\right.
\end{equation}
It is clear that $s_0A$ is doubly complementary. 
Applying Theorems~\ref{thm1} and \ref{thm2}, we have
\begin{proposition} \label{pro-cloaking} Let $d=2, 3$, $\delta > 0$,  and $D \subset \subset  B_{r_*} \setminus B_{r_2}$ and let $f \in L^2(\partial D)$. Assume that  $u_\delta$ and $v_\delta$ are defined by \eqref{def-ud} and \eqref{def-vd} where $A$ is given in \eqref{A-cloak}. 
There exists a sequence $\delta_n \to 0$ such that 
\begin{equation*}
\lim_{n \to \infty} E_{\delta_n} (u_{\delta_n}) = + \infty. 
\end{equation*}
Moreover, 
\begin{equation*}
v_{\delta_n} \to 0 \mbox{ weakly in } H^1(\Omega \setminus B_{r_3}). 
\end{equation*}
\end{proposition} 

\noindent {\bf Proof.} By Theorem~\ref{thm1} and Theorem~\ref{thm2}, it suffices to prove that there is no $W \in H^1(B_{r_*} \setminus B_{r_2})$ such that 
\begin{equation*}
\dive(A \nabla W) = f \mbox{ in } B_{r_*} \setminus B_{r_2} \quad \mbox{ and } \quad W= A \nabla W \cdot \eta = 0 \mbox{ on } \partial B_{r_2}. 
\end{equation*}
In fact, Theorems~\ref{thm1} and \ref{thm2} only deal with the case $f \in L^2(\Omega)$, however, the same results hold for $f$  stated here and  the proofs are unchanged. Suppose that this is not true, i.e., such a $W$ exists.
Since $\dive(A \nabla W ) = 0 $ in $(B_{r_*} \setminus B_{r_2}) \setminus \bar D$ and $W = A \nabla W \cdot \eta = 0$ on $\partial B_{r_2}$,  it follows from the unique continuation principle that 
$W = 0$ in $(B_{r_*} \setminus B_{r_2}) \setminus \bar D$. Hence $W = 0 $ in $D$ since $W \in H^1(B_{r_*} \setminus B_{r_2})$, $W = 0$ on $\partial D$,  and $\dive (A \nabla W) = 0$ in $D$. We deduce that $W = 0 $ in $ B_{r_*} \setminus B_{r_2}$. Hence $W = 0$ in $B_{r_*} \setminus B_{r_2}$.  This contradicts the fact that $\dive(A \nabla W) = f \neq 0$ in $B_{R_*} \setminus B_{r_2}$. \proofend

\appendix
\section{Appendix: Proof of Proposition~\ref{pro1}} \label{A}
\renewcommand{\theequation}{A\arabic{equation}}
\renewcommand{\thelemma}{A\arabic{lemma}}
  \setcounter{equation}{0}  
  \setcounter{lemma}{0}  

This appendix containing two subsections is devoted to  the proof of Proposition~\ref{pro1}. Some useful lemmas are established in the first section and  the proof of Propositions~\ref{pro1} is given in the second subsection. 

\subsection{Preliminaries}

In this section, we assume that 
\begin{equation*}
a \in [C^3(\overline{B_{R_2} \setminus B_{R_1}})]^{d \times d} \mbox{ is uniformly elliptic symmetric}, 
\end{equation*}
$$
K: B_{R_2}  \setminus B_{R_1} \to B_{R_3} \setminus B_{R_2} \mbox{ is defined by } K(x) = x R_2^2 / |x|^2, 
$$
and  $a_1$  is given by \eqref{def-a1}:
\begin{equation*}
a_1 = \left\{\begin{array}{cl} a &  \mbox{ in } B_{R_2} \setminus B_{R_1}, \\[6pt]
K_* a & \mbox{ in } B_{R_3} \setminus B_{R_2}, \\[6pt]
I & \mbox{ in } B_{R_1}. 
\end{array}\right.
\end{equation*}

We begin with 
\begin{lemma}\label{lem1}
Let $d=2, 3$, $v \in H^1(B_{R_3})$ be a solution to 
\begin{equation*}
\dive (a_1 \nabla v) = 0 \mbox{ in } B_{R_3}, 
\end{equation*}
and $w$ be the reflection of $v$ through $\partial B_{R_2}$ by $K^{-1}$, i.e., 
\begin{equation*}
w= v \circ K \mbox{ in } B_{R_2} \setminus B_{R_1}.
\end{equation*} 
Assume that 
\begin{equation}\label{pro-V}
v - w + c  = 0 \mbox{ on } \partial B_{R_1},
\end{equation}
for some $c \in \mC$. Then 
\begin{equation}\label{claim2-1-1-1}
v \mbox{ is constant and } c = 0.
\end{equation} 
\end{lemma}

\noindent{\bf Proof.} By considering the real part and the imaginary part separately, one may assume that $v$, $w$, and $c$ are real.  We first prove that $c = 0$. Assume that $c \neq 0$. From the definition of $w$ and \eqref{pro-V}, we have
\begin{equation}\label{max}
v(R_1 \sigma ) = v(R_3 \sigma) - c \quad \forall \, \sigma \in \partial B_1. 
\end{equation}
By the standard theory of elliptic equations,
\begin{equation*}
\sup_{\sigma \in \partial B_1 } |v(R_1 \sigma)| < + \infty, 
\end{equation*}
which implies, by \eqref{max},   
\begin{equation}\label{boundedness}
\sup_{\sigma \in \partial B_1 }|v(R_3 \sigma)| < + \infty. 
\end{equation}
Set, for $t \in \mR$,  
\begin{equation*}
b(t) = \sup_{\sigma \in \partial B_{1}} |v(R_3 \sigma) + t|. 
\end{equation*}
Applying the maximum principle,  we derive from \eqref{max} that
\begin{equation*}
\sup_{\sigma \in \partial B_1} |v(R_3 \sigma) + t| = \sup_{\sigma \in \partial B_1} |v(R_1 \sigma) + (t + c)| \le  \sup_{\sigma \in \partial B_1} |v(R_3 \sigma) + (t + c)|; 
\end{equation*}
this implies
\begin{equation*}
b(t) \le b(t + c). 
\end{equation*}
It follows that 
\begin{equation*}
b(-mc) \le b(0) \quad \forall \, m \ge 1:
\end{equation*}
we have a contradiction by \eqref{boundedness}. Hence $c = 0$.  From \eqref{max} and the maximum principle, we derive that $v$ is constant.  The proof is complete.   \proofend

\medskip

We also have 
\begin{lemma}\label{lem2}
Let $d = 2, 3$, $v \in H^1(B_{R_3})$ be a solution to 
\begin{equation*}
\dive (a_1 \nabla v) = 0 \mbox{ in } B_{R_3}, 
\end{equation*}
and $w$ be the reflection of $v$ through $\partial B_{R_2}$ by $K^{-1}$, i.e., $w = v \circ K$ in $B_{R_2} \setminus B_{R_1}$.  Set 
\begin{equation*}
V = v+ w. 
\end{equation*}
Assume that 
\begin{equation*}
a \nabla V \cdot \eta = c \mbox{ on } \partial B_{R_1}, 
\end{equation*}
for some $c \in \mC$. Then  
\begin{equation}\label{claim2-2}
v \mbox{ is constant and } c = 0.
\end{equation} 
\end{lemma}

\noindent{\bf Proof.} From the definition of $a_1$, by Lemma~\ref{lem-TO},  we have 
\begin{equation}\label{eq-V}
\dive (a \nabla V) = 0 \mbox{ in } B_{R_{2}} \setminus B_{R_1}
\end{equation}
and 
\begin{equation}\label{bd-V}
V = 2 v \quad \mbox{ and } \quad a \nabla V \cdot \eta = 0 \mbox{ on } \partial B_{R_2}. 
\end{equation}
Integrating \eqref{eq-V} in $B_{R_2} \setminus B_{R_1}$ and using \eqref{bd-V}, we obtain 
\begin{equation*}
\int_{\partial B_{R_1}} a \nabla V \cdot \eta = 0; 
\end{equation*}
which implies $c =0$. Hence, 
\begin{equation*}
a \nabla V \cdot \eta = 0  \mbox{ on } \partial B_{R_1} \cup \partial B_{R_2}. 
\end{equation*} 
It follows from \eqref{eq-V} that  $V$ is  constant in $ B_{R_2} \setminus B_{R_1}$. We derive from \eqref{bd-V} that  $v$ is constant on $\partial B_{R_2}$; hence $v$ is constant in $B_{R_3}$ by the unique continuation principle.  \proofend

\medskip 
The following lemma is one of the main ingredients in the proof of statement 1) of  Proposition~\ref{pro1} in two dimensions. 

\begin{lemma} \label{lem-density-2}  Let $d = 2$,  $v_{\ell, \pm} \subset H^1(B_{R_3})$ ($\ell \ge 1$) be the unique solution to  
\begin{equation}\label{def-vl}
\dive (a_1 \nabla v_{\ell, \pm} ) = 0 \mbox{ in } B_{R_3} \quad \mbox{ and } \quad v_{\ell, \pm} = e^{ \pm i \ell \theta} \mbox{ on } \partial B_{R_3},  
\end{equation}
and set 
\begin{equation*}
v_0 = 1 \mbox{ in } B_{R_3}. 
\end{equation*}
Define $w_{\ell, \pm} \in H^1(B_{R_2} \setminus B_{R_1})$ ($\ell \ge 1$)  the reflection of $v_{\ell, \pm}$ through $\partial B_{R_2}$ by $K^{-1}$, i.e., 
\begin{equation}\label{def-wl}
w_{\ell, \pm} = v_{\ell, \pm} \circ K \mbox{ in } B_{R_2} \setminus B_{R_1}, 
\end{equation}
and denote  $w_0 \in H^1(B_{R_2} \setminus B_{R_1})$  the unique solution to  
\begin{equation}\label{def-w0}
\dive (a_1 \nabla w_0) = 0 \mbox{ in } B_{R_2} \setminus B_{R_1}, \quad w_{0} = 1 \mbox{ on } \partial B_{R_2}, \quad \mbox{ and } \quad w_0 = 0 \mbox{ on } \partial B_{R_1}.   
\end{equation}
Then  
\begin{equation}\label{dense-1}
\{v_0 - w_0 \} \cup \Big\{v_{\ell, \pm} - w_{\ell, \pm}; \; \ell \ge 1 \Big\} \mbox{ is a dense subset of  } H^{1/2}(\partial B_{R_1}). 
\end{equation}
\end{lemma}

\noindent{\bf Proof.}  Let $G(x, y)$ be the fundamental solution to the equation $\dive (a_1 \nabla u) = 0$ in $B_{R_3}$ with respect to the zero Dirichlet boundary condition, i.e., 
\begin{equation*}
\dive_y (a_1 (y) \nabla_y G(x, y)) = \delta_x \mbox{ in } B_{R_3} \quad \mbox{ and } \quad G(x, y) = 0 \mbox{ on } \partial B_{R_3}. 
\end{equation*}
We have, by the Green formula, 
\begin{equation}\label{representation}
v_{\ell, \pm}(x) = \int_{\partial B_{R_3}} a_1(y)  \nabla_y G(x, y) \cdot \eta_y \,  v_{\ell, \pm}(y)\, dy 
\end{equation}
and, see e.g., \cite{DolzmannMuller} \footnote{The corresponding result in three dimensions can be found in \cite{GruterWidman}.},  
\begin{equation}\label{Green1}
|G(x, y)| \le C \mbox{ for } x \in B_{R_2}, y \in B_{R_3} \setminus B_{(R_2 + R_3)/2}. 
\end{equation}
Here and in what follows in this proof, $C$ denotes a positive constant independent of $x$, $y$, and $\ell$. 
It follows from \eqref{Green1} that, for $|\alpha| \le 2$, (see,  e.g.,  \cite[Theorems 6.2 and 6.6]{GilbargTrudinger})  
\begin{equation}\label{derivative-estimate}
|D^\alpha G(x, y)| \le C \mbox{ for } x \in B_{R_2}, y \in B_{R_3} \setminus B_{(R_2 + R_3)/2},  
\end{equation}
since $a_1 \in [C^3(\overline{B_{R_3} \setminus B_{(R_2 + R_3)/2}})]^{2 \times 2}$. 
A combination of \eqref{representation} and \eqref{derivative-estimate} yields 
\begin{equation}\label{est-derivative}
|\nabla v_{\ell, \pm} (x) | \le C/ \ell \quad \mbox{ for } x \in B_{R_3}, \; \ell \ge 1. 
\end{equation}
We claim that, for $\ell_0 \in \mN$ large enough, 
\begin{equation}\label{basis}
\{ e^{\pm i \ell \theta}; 0 \le \ell \le \ell_0 -1\} \cup \{v_{\ell, \pm} - w_{\ell, \pm}; \ell \ge \ell_0 \} \mbox{ is dense in } H^{1/2}(\partial B_{R_1}). 
\end{equation}
Consider the linear transformations
\begin{equation*}
{\cal J}, {\cal P}: H^{1/2}(\partial B_{R_1}) \to H^{1/2}(\partial B_{R_1})
\end{equation*}
defined as follows
\begin{equation*}
{\cal J} (e^{ \pm i \ell \theta}) = \left\{\begin{array}{cl} - e^{\pm i \ell \theta} & \mbox{ if } 0 \le \ell < \ell_0, \\[6pt]
v_{\ell, \pm} - w_{\ell, \pm} & \mbox{ if } \ell \ge \ell_0, 
\end{array}\right.
\end{equation*}
and
\begin{equation*}
{\cal P} (e^{ \pm i \ell \theta}) = \left\{\begin{array}{cl} 0 & \mbox{ if } 0 \le \ell < \ell_0, \\[6pt]
v_{\ell, \pm} & \mbox{ if } \ell \ge \ell_0. 
\end{array}\right.
\end{equation*}
Since $w_{\ell, \pm} = e^{\pm i \ell \theta}$ on $\partial B_{R_1}$, it follows that  
\begin{equation*}
{\cal J} = -  {\cal I} + {\cal P}, 
\end{equation*}
where ${\cal I}$ denotes the identity transformation.  

\medskip
Given $f \in H^{1/2}(\partial B_{R_1})$,   $f$ can be represented by 
\begin{equation*}
f = \alpha_0 + \sum_{\ell =1}^\infty \sum_{\pm} \alpha_{\ell, \pm} e^{\pm i \ell \theta} \mbox{ on } \partial B_{R_1}, 
\end{equation*}
for some $\alpha_0, \alpha_{\ell, \pm} \in \mC$ ($\ell \ge 1$).  We have 
\begin{equation*}
|\alpha_0|^2 +\sum_{\ell \ge 1} \sum_{\pm} \ell |\alpha_{\ell, \pm}|^2 \le C \| f\|_{H^{1/2}(\partial B_{R_1})}^2. 
\end{equation*}
From the definition of ${\cal P}$, 
\begin{equation*}
{\cal P}(f) = \sum_{\ell \ge \ell_0} \sum_{\pm} \alpha_{\ell, \pm}  v_{\ell, \pm} \mbox{ on } \partial B_{R_1}. 
\end{equation*}
We derive from \eqref{est-derivative} that
\begin{align*}
\| {\cal P}(f) \|_{H^{1/2}(\partial B_{R_1})} \le C \sum_{\ell \ge \ell_0} \sum_{\pm} |\alpha_{\ell, \pm}|/ \ell & \le C \left(\sum_{\ell \ge \ell_0} \sum_{\pm}  \ell|\alpha_{\ell, \pm}|^2 \right)^{1/2} \left(\sum_{\ell \ge \ell_0}  \sum_{\pm} 1/\ell^{3} \right)^{1/2} \\[6pt]
& \le C \ell_0^{-1} \| f\|_{H^{1/2}}. 
\end{align*}
Thus,  for $\ell_0$ large enough, $\|{\cal P}\| \le 1/2$. Hence ${\cal J}$ is invertible and  \eqref{basis} follows.

\medskip
Fix $\ell_0$ such that \eqref{basis} holds. Using \eqref{basis}, we derive   that the dimension of  the orthogonal complement of $ \{v_{\ell, \pm} - w_{\ell, \pm}; \ell \ge \ell_0 \} $ in $H^{1/2} (\partial B_{R_1})$ is less than or equal to $2 \ell_0 - 1$. Hence, to obtain the conclusion, it suffices to prove that 
\begin{equation}\label{claim1}
\{U_0 \} \cup \big\{ U_{\ell, \pm} \big\}_{1 \le \ell < \ell_0} \mbox{ is linearly independent in } H^{1/2}(\partial B_{R_1}), 
\end{equation}
where $U_0$ and $U_{\ell, \pm}$ ($1 \le  \ell < \ell_0$) are respectively the projection of $v_0 - w_0$ and $v_{\ell, \pm} - w_{\ell, \pm}$ into $\Big(\mbox{span}\{v_{\ell, \pm} - w_{\ell, \pm}; \ell \ge \ell_0 \}\Big)^\perp$ with respect to $H^{1/2}(\partial B_{R_1})$ scalar product.  Indeed, let  $\alpha_0, \alpha_{\ell, \pm} \in \mC$ ($1 \le \ell  < \ell_0$) be such that 
\begin{equation}\label{assumption1}
\alpha_0 U_0 + \sum_{\ell = 1}^{\ell_0 -1 } \sum_{\pm} \alpha_{\ell, \pm} U_{\ell, \pm} = 0 \mbox{ on } \partial B_{R_1}. 
\end{equation}
We prove that $\alpha_0 = \alpha_{\ell, \pm} = 0$ for $1 \le \ell \le \ell_0 -1$. From \eqref{assumption1}, we have 
\begin{equation*}
\alpha_0 (v_0 - w_0) + \sum_{\ell = 1}^{\ell_0 - 1}  \sum_{\pm} \alpha_{\ell, \pm} (v_{\ell, \pm} - w_{\ell, \pm}) = v - w \mbox{ on } \partial B_{R_1}, 
\end{equation*}
for some $v \in \mbox{closure}\Big\{ \mbox{span}\{v_{\ell, \pm}; \ell \ge \ell_0 \}\Big\}$ with respect to $H^1(B_{R_3})$-norm. Here $w$ is the reflection of $v$ through $\partial B_{R_2}$ by $K^{-1}$. Set 
\begin{equation}\label{def-V}
V = \sum_{\ell = 1}^{\ell_0 - 1} \sum_{\pm} \alpha_{\ell, \pm} v_{\ell, \pm} - v \mbox{ in } B_{R_3}, 
\end{equation}
and denote $W$ the reflection of $V$ through $\partial B_{R_2}$ by $K^{-1}$. It follows that 
\begin{equation*}
\alpha_0 (v_0 - w_0) + V - W = 0 \mbox{ on } \partial B_{R_1}. 
\end{equation*}
Applying Lemma~\ref{lem1}, we have 
\begin{equation*}
\alpha_0 = 0  \quad \mbox{ and } \quad V \mbox{ is  constant}. 
\end{equation*}
We derive from the definition of $V$ in \eqref{def-V} that 
\begin{equation*}
\alpha_{\ell, \pm} = 0 \mbox{ for } 1 \le \ell \le \ell_0 - 1. 
\end{equation*}
The proof of \eqref{claim1} is complete. \proofend

\medskip
For $D$ an open subset $D$ of $\mR^d$, we denote 
\begin{equation*}
H^1_{\sharp}(D) = \Big\{v \in H^1(D); \; \int_{D} v = 0 \Big\}. 
\end{equation*}
The following result, which  is a variant of Lemma~\ref{lem-density-2} when the Neumann data on $\partial B_{R_1}$ is considered,  plays an important role in the proof of statement 2) of Proposition~\ref{pro1}. 

\begin{lemma} \label{lem-density-2-1}  Let $d = 2$ and let  $v_{\ell, \pm} \subset H^1_{\sharp}(B_{R_3})$ ($\ell \ge 1$) be the unique solution to  
\begin{equation}\label{def-vl-1}
\dive (a_1 \nabla v_{\ell, \pm} ) = 0 \mbox{ in } B_{R_3} \quad \mbox{ and } \quad a_1 \nabla v_{\ell, \pm} \cdot \eta = e^{ \pm i \ell \theta} \mbox{ on } \partial B_{R_3},  
\end{equation}
Define $w_{\ell, \pm} \in H^1_{\sharp}(B_{R_2} \setminus B_{R_1})$  the reflection of $v_{\ell, \pm}$ through $\partial B_{R_2}$ by $K^{-1}$, i.e., 
\begin{equation}\label{def-wl-1}
w_{\ell, \pm} = v_{\ell, \pm} \circ K \mbox{ in } B_{R_2} \setminus B_{R_1}. 
\end{equation}
We have
\begin{equation}\label{dense-2}
\Big\{ 1 \Big\} \cup \Big\{a \nabla (v_{\ell, \pm} + w_{\ell, \pm}) \cdot \eta; \; \ell \ge 1 \Big\} \mbox{ is a dense subset of } H^{-1/2}(\partial B_{R_1}). 
\end{equation}
\end{lemma}

\begin{remark} \fontfamily{m} \selectfont
Since $\dsp \int_{\partial B_{R_3}} e^{\pm i \ell \theta} = 0$ for $\ell \ge 1$, it follows that  $v_{\ell, \pm}$ is well-defined. 
\end{remark}

\medskip
\noindent {\bf Proof.} The proof of Lemma~\ref{lem-density-2-1} is in the same spirit of the one of Lemma~\ref{lem-density-2}.  As in the proof of Lemma~\ref{lem-density-2}, we also reach 
\begin{equation}
\{1\} \cup \{ e^{\pm i\ell \theta}; 1 \le \ell < \ell_0 \} \cup \{ a \nabla (v_\ell + w_\ell ) \cdot \eta; \; \ell \ge \ell_0 \} \mbox{ is dense in } H^{-1/2}(\partial B_{R_1}), 
\end{equation}
for some $\ell_0 > 1$ (large). It follows that the dimension of  the orthogonal complement of \\
$\mbox{closure} \Big\{\mbox{span} \{ a \nabla (v_\ell + w_\ell) \cdot \eta ; \ell \ge \ell_0 \} \Big\}$ in $H^{-1/2} (\partial B_{R_1})$ is less than or equal to $2\ell_0 - 1$. Hence, to obtain the conclusion, it suffices to prove that 
\begin{equation}\label{claim2-1}
\{U_0\} \cup \big\{ U_{\ell, \pm} \big\}_{1 \le \ell < \ell_0} \mbox{ is independent in } H^{-1/2}(\partial B_{R_1}), 
\end{equation}
where  $U_0 = 1$ and $U_{\ell, \pm}$ ($1 \le   \ell < \ell_0$) is the projection of $a \nabla (v_{\ell, \pm} + w_{\ell, \pm}) \cdot \eta $ into \\
 $\Big(\mbox{closure} \Big\{\mbox{span}\{ a\nabla (v_{\ell, \pm} + w_{\ell, \pm}) \cdot \eta ; \ell \ge \ell_0 \} \Big\}\Big)^\perp$ with respect to  $H^{-1/2}(\partial B_{R_1})$ scalar product.

\medskip
 Let  $\alpha_0, \alpha_{\ell, \pm} \in \mC$ ($1 \le \ell \le \ell_0 -1$) be such that 
\begin{equation}\label{assumption-2-1}
\alpha_0 + \sum_{\ell = 1}^{\ell_0 -1 } \sum_{\pm}\alpha_{\ell, \pm} U_{\ell, \pm} = 0 \mbox{ on } \partial B_{R_1}. 
\end{equation}
We prove that $\alpha_0 = \alpha_{\ell, \pm} = 0$ for $1 \le \ell \le \ell_0 -1$. From \eqref{assumption-2-1}, we have 
\begin{equation}\label{cl2-1}
\alpha_0 + \sum_{\ell = 1}^{\ell_0 - 1} \sum_{\pm} \alpha_{\ell, \pm} a \nabla (v_{\ell, \pm} + w_{\ell, \pm}) \cdot \eta = a \nabla (v + w) \cdot \eta \mbox{ on } \partial B_{R_1}, 
\end{equation}
for some $v \in \mbox{closure}\Big\{ \mbox{span}\{v_{\ell, \pm}; \ell \ge \ell_0 \}\Big\}$ in $H^1_{\sharp}(B_{R_3})$. Here $w$ is the reflection of $v$ through $\partial B_{R_2}$ by $K^{-1}$. Set 
\begin{equation}\label{def-V-2}
V = \sum_{\ell = 1}^{\ell_0 - 1} \sum_{\pm} \alpha_{\ell, \pm} v_{\ell, \pm} - v \mbox{ in } B_{R_3}, 
\end{equation}
and denote $W$ the reflection of $V$ through $\partial B_{R_2}$ by $K^{-1}$. It follows from \eqref{cl2-1} that 
\begin{equation*}
\alpha_0 +  a \nabla( V + W) \cdot \eta = 0 \mbox{ on } \partial B_{R_1}. 
\end{equation*}
Applying Lemma~\ref{lem2}, we have 
\begin{equation*}
\alpha_0 = 0  \quad \mbox{ and } \quad V \mbox{ is  constant}. 
\end{equation*}
Hence $V=0$ since $V \in H^1_{\sharp}(B_{R_3})$. We derive from the definition of $V$ in \eqref{def-V-2} and of $v_{\ell, \pm}$ that 
\begin{equation*}
\alpha_{\ell, \pm} = 0 \mbox{ for } 1 \le \ell \le \ell_0 - 1. 
\end{equation*}
The proof of \eqref{claim2-1} is complete.  \proofend

%

\medskip 
Here are variants of Lemmas~\ref{lem-density-2} and  \ref{lem-density-2-1} in three dimensions. The first one is the variant of Lemma~\ref{lem-density-2}. 

\begin{lemma} \label{lem-density-3}  Let $d = 3$ and  let  $v^{k}_{\ell} \subset H^1(B_{R_3})$ ($\ell \ge 1, - \ell \le k \le \ell$) be the unique solution to  
\begin{equation}\label{def-vl-4}
\dive (a_1 \nabla v^k_{\ell} ) = 0 \mbox{ in } B_{R_3} \quad \mbox{ and } \quad v^k_{\ell} = Y^k_\ell \mbox{ on } \partial B_{R_3}. 
\end{equation}
and set $v^0_0 = 1$.  Here $Y^k_\ell$ is the spherical harmonic function of degree $\ell$ and of order $k$. 
Define $w^k_\ell \in H^1(B_{R_2} \setminus B_{R_1})$  the reflection of $v^k_\ell$ through $\partial B_{R_2}$ by $K^{-1}$, i.e., 
\begin{equation}\label{def-wl-4}
w^k_\ell = v^k_\ell \circ K \mbox{ in } B_{R_2} \setminus B_{R_1}, 
\end{equation}
and denote  $w_0^0 \in H^1(B_{R_2} \setminus B_{R_1})$  the unique solution to  
\begin{equation}\label{def-w0-4}
\dive (a_1 \nabla w^0_0) = 0 \mbox{ in } B_{R_2} \setminus B_{R_1}, \quad w^0_{0} = 1 \mbox{ on } \partial B_{R_2}, \quad \mbox{ and } \quad w^0_0 = 0 \mbox{ on } \partial B_{R_1}.   
\end{equation}
We have
\begin{equation}\label{dense-1-4}
\Big\{v^k_\ell - w^k_\ell; \; \ell \ge 0, - \ell \le k \le \ell \Big\} \mbox{ is a dense subset of  } H^{1/2}(\partial B_{R_1}). 
\end{equation}
\end{lemma}

\noindent{\bf Proof.} The proof is similar to the one of Lemma~\ref{lem-density-2}. The details are left to the reader. \proofend

\medskip 
The second one is the variant of Lemma~\ref{lem-density-2-1}. 

\begin{lemma} \label{lem-density-3-1}  Let $d = 3$ and  let  $v_{k, \ell} \subset H^1_{\sharp}(B_{R_3})$ ($\ell \ge 1, - \ell \le k \le \ell$) be the unique solution to  
\begin{equation}\label{def-vl-5}
\dive (a_1 \nabla v^k_{\ell} ) = 0 \mbox{ in } B_{R_3} \quad \mbox{ and } \quad a_1 \nabla v^k_{\ell} \cdot \eta = Y^k_\ell \mbox{ on } \partial B_{R_3}. 
\end{equation}
Define $w^k_\ell \in H^1_{\sharp}(B_{R_2} \setminus B_{R_1})$ ($\ell \ge 1$)  the reflection of $v^k_\ell$ through $\partial B_{R_2}$ by $K^{-1}$, i.e., 
\begin{equation}\label{def-wl-5}
w^k_\ell = v^k_\ell \circ K \mbox{ in } B_{R_2} \setminus B_{R_1}, 
\end{equation}
 We have
\begin{equation}\label{dense-2-5}
\{1\} \cup \Big\{a_1 \nabla (v^k_\ell +  w^k_\ell) \cdot \eta ; \; \ell \ge 1, - \ell \le k \le \ell \Big\} \mbox{ is a dense subset of  } H^{-1/2}(\partial B_{R_1}). 
\end{equation}
\end{lemma}

\medskip

\noindent{\bf Proof.} Since $\int_{\partial B_{R_3}} Y^m_{\ell} = 0$ for $\ell \ge 1$ and $-\ell \le k \le \ell $, it follows that  $v^k_{\ell}$ is well-defined. 
The proof is similar to the one of Lemma~\ref{lem-density-2-1}. The details are left to the reader. \proofend

%

\subsection{Proof of Proposition~\ref{pro1}.}

Statements \eqref{pro-wl-1} and \eqref{pro-wl-2} are consequences of Lemma~\ref{lem-TO}. It remains to prove statements 1), 2), and 3). The proof is now divided into two steps. 

\medskip
\noindent{\bf Step 1:}  We prove that if one of  statements 1), 2), and 3) of  Proposition~\ref{pro1} hold for a (particular) dense set $(v_\ell)_{\ell \ge 0}$, then it also  holds for all dense sets $(v_\ell)_{\ell \ge 0}$.

We will only discuss this fact for statement $1)$, the  other cases follows similarly.  Assume that statement $1)$ holds for a specific sequence of $\{v_\ell \}_{\ell  \ge 0}$ which satisfies the assumptions of Proposition~\ref{pro1}. 
We will prove that statement $1)$ holds for any sequence $\{ \hat v_\ell \}_{\ell  \ge 0}$ satisfying the assumptions of Proposition~\ref{pro1}. Let 
$v \in H^1(B_{R_2} \setminus B_{R_1})$ be such that $\dive (a \nabla v) = 0$ in $B_{R_2} \setminus B_{R_1}$ and $v = 0$ on $\partial B_{R_2}$. 
For $\eps > 0$, there exist $\ell_\eps > 0$ and $(\alpha_{\ell})_{0}^{\ell_\eps} \subset \mC$ such that 
\begin{equation}\label{bb-0}
\|v - \sum_{0}^{\ell_\eps} \alpha_\ell (v_\ell - w_\ell) \|_{H^1(B_{R_2} \setminus B_{R_1})} \le \eps. 
\end{equation}
since statement 1 holds for $(v_\ell)$. On the other hand, there exist $\hat \ell_\eps$  and $(\hat \alpha_\ell)_0^{\hat \ell_\eps} \subset \mC$ such that 
\begin{equation*}
\| \sum_{0}^{\ell_\eps} \alpha_\ell v_\ell  - \sum_{0}^{\hat \ell_\eps} \hat \alpha_\ell \hat v_\ell \|_{H^{1/2}(\partial B_{R_3})} \le \eps, 
\end{equation*}
by the dense property of $\{\hat v_\ell\}_0^\infty$. This implies 
\begin{equation}\label{bb-1}
\| \sum_{0}^{\ell_\eps} \alpha_\ell v_\ell  - \sum_{0}^{\hat \ell_\eps} \hat \alpha_\ell \hat v_\ell \|_{H^1(B_{R_3})} \le \eps.  
\end{equation}
Let $\hat w_\ell$ be the reflection of $\hat v_\ell$ through $\partial B_{R_2}$ by $K^{-1}$ for $\ell \ge 1$. Note that if $w$ is the reflection of $v$ through  $\partial B_{R_2}$ by $K^{-1}$, then
\begin{equation}\label{bb-1-1}
\| w\|_{H^1(B_{R_2} \setminus B_{R_1})} \le C \| v\|_{H^1(B_{R_3})}. 
\end{equation}
Here and in what follows $C$ denotes a positive constant depending only on $a$, $R_1$, and $R_2$. 
A combination of \eqref{bb-1} and \eqref{bb-1-1} yields 
\begin{equation}\label{bb-2}
\| \sum_{1}^{\ell_\eps} \alpha_\ell w_\ell  - \sum_{1}^{\hat \ell_\eps} \hat \alpha_\ell \hat w_\ell + (\alpha_0 - \hat\alpha_0) \|_{H^1(B_{R_2} \setminus B_{R_1})} \le   C \eps. 
\end{equation}
We derive from  \eqref{bb-1} and \eqref{bb-2} that 
\begin{equation}\label{bb-3}
\| \sum_{1}^{\ell_\eps} \alpha_\ell (v_\ell - w_\ell)  - \sum_{1}^{\hat \ell_\eps} \hat \alpha_\ell (\hat v_\ell -  \hat w_\ell)  \|_{H^1(B_{R_2} \setminus B_{R_1})} \le  C \eps. 
\end{equation}
From \eqref{bb-0} and \eqref{bb-3}, we obtain 
\begin{equation*}
\|v - \sum_{1}^{\hat \ell_\eps} \hat \alpha_\ell (\hat v_\ell - \hat w_\ell) - \alpha_0(v_0 - w_0) \|_{H^1(B_{R_2} \setminus B_{R_1})} \le C \eps. 
\end{equation*}
Hence statement 1) holds for $(\hat v_\ell)$.

\medskip 
\noindent{\bf Step 2:} Proof of statements 1), 2), and 3).  

We only establish these statements in two dimensions. The three dimensional case follows similarly.  However, instead of applying Lemmas~\ref{lem-density-2} and \ref{lem-density-2-1}, one uses Lemmas~\ref{lem-density-3} and \ref{lem-density-3-1}.

\medskip 
Assume $d=2$.  Let  $v_{\ell, \pm} \subset H^1(B_{R_3})$ ($\ell \ge 1$) be the unique solution to  
\begin{equation}\label{def-vl-7}
\dive (a_1 \nabla v_{\ell, \pm} ) = 0 \mbox{ in } B_{R_3} \quad \mbox{ and } \quad v_{\ell, \pm} = e^{\pm i \ell \theta} \mbox{ on } \partial B_{R_3}, 
\end{equation}
and set 
\begin{equation}\label{v0}
v_0 = 1 \mbox{ in } B_{R_3}. 
\end{equation}
Let $w_{\ell, \pm} \in H^1(B_{R_2} \setminus B_{R_1})$ ($\ell \ge 1$)  be the reflection of $v_{\ell, \pm}$ through $\partial B_{R_2}$ by $K^{-1}$, i.e., 
\begin{equation}\label{def-wl-7}
w_{\ell, \pm} = v_{\ell, \pm} \circ K \mbox{ in } B_{R_2} \setminus B_{R_1}, 
\end{equation}
and denote  $w_0 \in H^1(B_{R_3} \setminus B_{R_2})$  the unique solution to  
\begin{equation*}
\dive (a \nabla w_0) = 0 \mbox{ in } B_{R_2} \setminus B_{R_1}, \quad w_{0} = 1 \mbox{ on } \partial B_{R_2}, \quad \mbox{ and } \quad w_0 = 0 \mbox{ on } \partial B_{R_1}.   
\end{equation*}
By Step 1, it suffices to prove the statements 1), 2), and 3) for  $\{v_0, w_0 \} \cup \{v_{\ell, \pm}, w_{\ell, \pm}\}_{\ell \ge 1}$.    

\medskip
\noindent{\bf Proof of statement 1).} This statement is a consequence of the fact that 
 $v = 0$ if  $v \in H^1(B_{R_2} \setminus B_{R_1})$ satisfies 
 \begin{equation}\label{orthogonal-1-0}
\dive(a \nabla v) = 0 \mbox{ in } B_{R_2} \setminus B_{R_1},  \quad v = 0 \mbox{ on } \partial B_{R_2},
\end{equation}
\begin{equation}\label{orthogonal-1}
\int_{B_{R_2} \setminus B_{R_1}} a \nabla v \nabla (\bar v_{\ell, \pm} - \bar w_{\ell, \pm}) = 0  \quad \forall \,  \ell \ge 1,  
\end{equation}
and 
\begin{equation}\label{orthogonal-1-1}
\int_{B_{R_2} \setminus B_{R_1}} a \nabla v \nabla (\bar v_0 - \bar w_0) = 0. 
\end{equation}
Indeed, using \eqref{orthogonal-1-0}, we derive from \eqref{orthogonal-1} and \eqref{orthogonal-1-1} that 
\begin{equation}\label{orthogonal-1-2}
\int_{\partial B_{R_1}} a \nabla v \cdot \eta \; (\bar v_{\ell, \pm} - \bar w_{\ell, \pm})  = 0 \quad \forall \, \ell \ge 1 
\end{equation}
and 
\begin{equation}\label{orthogonal-1-2-1}
\int_{\partial B_{R_1}} a \nabla v \cdot \eta \; (\bar v_0 - \bar w_0)  = 0. 
\end{equation}
Since, by Lemma~\ref{lem-density-2}, 
\begin{equation*}
\{v_0 - w_0 \} \cup \Big\{v_{\ell, \pm} - w_{\ell, \pm}; \; \ell \ge 1 \Big\} \mbox{ is dense in  } H^{1/2}(\partial B_{R_1}). 
\end{equation*}
 it follows from \eqref{orthogonal-1-2} and  \eqref{orthogonal-1-2-1}  that \begin{equation*}
a \nabla v \cdot \eta = 0  \mbox{ on } \partial B_{R_1}. 
\end{equation*} 
We derive from \eqref{orthogonal-1-0}  that $v =0$ in $B_{R_2} \setminus B_{R_1}$: statement 1) is proved.

\medskip
\noindent {\bf Proof of statement 2):} 
This statement is a consequence of the fact that 
 $v$ is constant  if  $v \in H^1(B_{R_2} \setminus B_{R_1})$ satisfies
\begin{equation}\label{orthogonal-2-0}
\dive(a \nabla v) = 0 \mbox{ in } B_{R_2} \setminus B_{R_1}, \quad a \nabla v \cdot \eta = 0 \mbox{ on } \partial B_{R_2}, 
\end{equation}
and 
\begin{equation}\label{orthogonal-2}
\int_{B_{R_2} \setminus B_{R_1}} a \nabla v \nabla (\bar v_{\ell, \pm} + \bar w_{\ell, \pm})= 0  \quad \forall \,  \ell \ge 1. 
\end{equation}
Indeed, since $a \nabla v_{\ell, \pm} \cdot \eta = - a \nabla w_{\ell, \pm} \cdot \eta$ on $\partial B_{R_2}$ for $\ell \ge 1$ \eqref{pro-wl-2},  it follows from \eqref{orthogonal-2} that 
\begin{equation}\label{orthogonal-2-1}
\int_{\partial B_{R_1}}  a \nabla (\bar v_{\ell, \pm} + \bar w_{\ell, \pm})  \cdot \eta \; v = 0 \quad \forall \, \ell \ge 1. 
\end{equation}
By  Lemma~\ref{lem-density-2-1} and Step 1, 
\begin{equation}\label{hoho}
\Big\{ 1 \Big\} \cup \Big\{a \nabla (v_{\ell, \pm} + w_{\ell, \pm}) \cdot \eta; \; \ell \ge 1 \Big\} \mbox{ is a dense subset of } H^{-1/2}(\partial B_{R_1}). 
\end{equation}
We derive from \eqref{orthogonal-2-1} that 
\begin{equation*}
 v  \mbox{ is constant on } \partial B_{R_1},   
\end{equation*} 
This implies, by \eqref{orthogonal-2-0},  
\begin{equation*}
v \mbox{ is constant in } B_{R_2} \setminus B_{R_1}. 
\end{equation*}
Statement 2) is proved.

\medskip
\noindent {\bf Proof of statement 3):}  
This statement is a consequence of the fact that 
 $v$ is constant  if  $v \in H^1(B_{R_2} \setminus B_{R_1})$ satisfies
\begin{equation}\label{orthogonal-0}
\dive(a \nabla v) = 0 \mbox{ in } B_{R_2} \setminus B_{R_1},
\end{equation}
\begin{equation}\label{orthogonal-11}
\int_{B_{R_2} \setminus B_{R_1}} a \nabla v \nabla \bar v_{\ell, \pm} = \int_{B_{R_2} \setminus B_{R_1}} a \nabla v \nabla \bar w_{\ell, \pm} = 0  \quad \forall \,  \ell \ge 1,  
\end{equation}
and 
\begin{equation}\label{orthogonal-11-1}
 \int_{B_{R_2} \setminus B_{R_1}} a \nabla v \nabla \bar v_0 = \int_{B_{R_2} \setminus B_{R_1}} a \nabla v \nabla \bar w_0 = 0. 
\end{equation}
In fact, a combination of \eqref{orthogonal-0}, \eqref{orthogonal-11}, and  \eqref{orthogonal-11-1} yields
\begin{equation}\label{tt1-pro}
\int_{\partial B_{R_2} \cup \partial B_{R_1}} a \nabla v \cdot \eta \; \bar  v_{\ell, \pm} = \int_{\partial B_{R_2} \cup \partial B_{R_1}} a \nabla v \cdot \eta \;  \bar w_{\ell, \pm}  = 0 \quad \forall \, \ell \ge 1 
\end{equation}
and 
\begin{equation}\label{tt2-pro}
\int_{\partial B_{R_2} \cup \partial B_{R_1}} a \nabla v \cdot \eta \;  \bar v_0  = \int_{\partial B_{R_2} \cup \partial B_{R_1}} a \nabla v \cdot \eta \;  \bar w_0  = 0.  
\end{equation}
Since $v_0 = w_0 =1$ and  $v_{\ell, \pm} = w_{\ell, \pm}$ on $\partial B_{R_2}$ for $\ell \ge 1$, it follows from \eqref{tt1-pro} that 
\begin{equation}\label{orthogonal-3}
\int_{\partial B_{R_1}} a \nabla v \cdot \eta \; (\bar v_{\ell, \pm} - \bar w_{\ell, \pm})  = 0 \quad \forall \, \ell \ge 1, 
\end{equation}
and, since $w_0 = 0$ on $\partial B_{R_1}$,  
\begin{equation}\label{orthogonal-3-1}
\int_{\partial B_{R_1}} a \nabla v \cdot \eta  = 0. 
\end{equation}
From \eqref{pro-wl-1}, \eqref{orthogonal-11},  and the symmetry of $a$, we also have
\begin{equation*}
\int_{\partial B_{R_2} \cup \partial B_{R_1}} a \nabla \bar v_{\ell, \pm} \cdot \eta \; \bar  v = \int_{\partial B_{R_2} \cup \partial B_{R_1}} a \nabla \bar w_{\ell, \pm} \cdot \eta \;  v  = 0 \quad \forall \, \ell \ge 1;  
\end{equation*}
which yields, since $a \nabla v_{\ell, \pm} \cdot \eta = - a \nabla w_{\ell, \pm} \cdot \eta$ for $\ell \ge 1$, 
\begin{equation}\label{orthogonal-4}
\int_{\partial B_{R_1}}  a \nabla (\bar v_{\ell, \pm} + \bar w_{\ell, \pm})  \cdot \eta \; v = 0 \quad \forall \, \ell \ge 1. 
\end{equation}
Using Lemma~\ref{lem-density-2} and \eqref{hoho}, we derive  from \eqref{orthogonal-3}, \eqref{orthogonal-3-1},   and \eqref{orthogonal-4} that  
\begin{equation}\label{orthogonal-5}
a \nabla v \cdot \eta = 0 \quad \mbox{ and } \quad v - \mint_{\partial B_{r_1}} v = 0 \mbox{ on } \partial B_{R_1}. 
\end{equation} 
A combination of \eqref{orthogonal-0} and \eqref{orthogonal-5} yields $v$ is constant  in $B_{R_2} \setminus B_{R_1}$ by the unique continuation principle.  Statement 3) is proved. 
The proof is complete. \proofend

\begin{remark}  \fontfamily{m} \selectfont In Proposition~\ref{pro1}, if one assumes in addition  that $\{v_\ell\}_0^\infty$ is a basis of $H^{1/2}(\partial B_{R_3})$. Then 
\begin{equation*}
\left|\begin{array}{c}
\{ v_\ell, w_\ell; \; \ell \ge 0 \} \mbox{ is finitely linearly independent in } H^1(B_{R_2} \setminus B_{R_1}), \\[6pt]
\{ v_\ell; \; \ell \ge 0 \} \mbox{ is finitely linearly independent in } H^{1/2}(\partial B_{R_2}), \\[6pt]
\{1 \} \cup \{ a \nabla w_\ell \cdot \eta; \; \ell \ge 1 \} \mbox{ is finitely linearly independent in } H^{-1/2}(\partial B_{R_2}). 
\end{array}\right.
\end{equation*}
These facts can be derived from Lemma~\ref{lem1}. 
\end{remark}

\bigskip

\noindent{\bf Acknowledgment: } The author thanks Graeme Milton for interesting discussions 
on the subject. The author also thanks Boris Buffoni for  useful discussions
which help him to simplify several points in the proof of Lemma~\ref{lem-density-2}. The author thanks the referees for helpful comments. 
The research was partially supported by NSF grant DMS-1201370 and by the Alfred P. Sloan Foundation.

\bibliographystyle{amsplain}

\providecommand{\bysame}{\leavevmode\hbox to3em{\hrulefill}\thinspace}
\providecommand{\MR}{\relax\ifhmode\unskip\space\fi MR }
\providecommand{\MRhref}[2]{%
  \href{http://www.ams.org/mathscinet-getitem?mr=#1}{#2}
}
\providecommand{\href}[2]{#2}

\end{document}